\documentclass[11pt]{article}
 \usepackage{benstyle}
 
\usepackage[font={small}]{caption}

\title{Infinite-dimensional compressed sensing and function interpolation}
\author{Ben Adcock \\ Department of Mathematics \\ Simon Fraser University \\ Canada}

\begin{document}

\maketitle

\begin{abstract}
We introduce and analyze a framework for function interpolation using compressed sensing.  This framework -- which is based on weighted $\ell^1$ minimization -- does not require \textit{a priori} bounds on the expansion tail in either its implementation or its theoretical guarantees, and in the absence of noise leads to genuinely interpolatory approximations.  We also establish a new recovery guarantee for compressed sensing with weighted $\ell^1$ minimization based on this framework.  This guarantee has several key benefits.  First, unlike existing results, it is sharp (up to constants and log factors) for large classes of functions regardless of the choice of weights.  Second, by examining the measurement condition in the recovery guarantee, we are able to suggest a good overall strategy for selecting the weights.  In particular, when applied to the important case of multivariate approximation with orthogonal polynomials, this weighting strategy leads to provably optimal estimates on the number of measurements required, whenever the support set of the significant coefficients is a so-called lower set.  Finally, this guarantee can also be used to theoretically confirm the benefits of alternative weighting strategies where the weights are chosen based on prior support information.  This provides a theoretical basis for a number of recent numerical studies showing the effectiveness of such approaches.
\end{abstract}

\section{Introduction}
Many problems in science and engineering require the approximation of smooth, multivariate functions from finitely-many pointwise samples.  Although a classical problem of approximation theory, recently there has been a renewed focus in this area, driven in part by applications in uncertainty quantification.  Problems in this area are typically high dimensional and place severe limitations on the number of measurements that can be acquired.  At the same time, developments in the field of compressed sensing (CS) have shown that it is often possible to recover high-dimensional vectors possessing certain low-dimensional structures from substantially reduced sets of linear measurements \cite{CandesRombergTao,donohoCS}.  It is known that smooth, high-dimensional functions have approximately sparse expansions in certain orthogonal systems (e.g.\ tensor Chebyshev or Legendre polynomials).  
Hence, in recent years there has been an increasing focus on applying the theory and techniques of CS to accurately compute such expansions
\cite{DoostanOwhadiSparse,HamptonDoostanCSPCE,MathelinGallivanCSPDErandom,NarayanZhouCCP,PengHamptonDoostantweighted,RauhutWardSpher,Rauhut,RauhutWardWeighted,YanGuoXui_l1UQ,KarniadakisUQCS}.

However, the application of CS to function approximation raises several issues.  First, standard CS concerns itself primarily with the recovery of sparse vectors in finite-dimensional vector spaces.  Functions, on the other hand, live in infinite-dimensional spaces.  Whilst they may be well approximated by finite sums in certain orthogonal polynomial systems, their expansion in a such system is typically infinite.  As one might expect, this mismatch presents a number of key practical and theoretical issues.  Second, the coefficients of a function in a polynomial basis are not just sparse, but tend to possess additional structure.  This raises questions about how to exploit such structure in the reconstruction process, and the expected theoretical benefits of doing so.  Third, it is well-known that na\"ive approximations of high-dimensional functions suffer from the \textit{curse of dimensionality}.  Therefore, a question of singular interest is whether, and the extent to which, CS provides a means to avoid this crucial issue.
 
In this paper we introduce a framework and series of recovery guarantees for the use of CS in approximating multivariate functions from limited numbers of pointwise samples.  Unlike most existing approaches to this problem -- which can loosely be described as discretizing first -- in our framework the recovery problem is first formulated as a (weighted) $\ell^1$ minimization problem in an infinite-dimensional space, and \textit{then} discretized.  This brings a number of benefits.  First, our techniques do not require \textit{a priori} estimates for the expansion tail, as is common in other approaches.  Second, in the absence of noise, our techniques lead to exactly interpolating approximations; a desirable property in general as well as for certain applications.  Third, much like our techniques, our theoretical results do not assume \textit{a priori} estimates for the expansion tail, as has been necessary in previous works.  Since infinite expansions are handled faithfully, continuing a line of investigation initiated in \cite{BAACHGSCS}, we refer to this framework as \textit{infinite-dimensional} compressed sensing.

In order to exploit the structure of polynomial coefficients of smooth functions we use a weighted $\ell^1$ minimization approach.  Incorporating weights into the reconstruction problem has received some attention of late, especially for polynomial approximations, due to its potential for enhancing accuracy; see \cite{Adcockl1pointwise,PengHamptonDoostantweighted,KarniadakisUQCS}, as well as \cite{RauhutWardWeighted} for some theoretical analysis.  In this paper, we introduce a new recovery guarantee for CS for very general choices of optimization weights.  For the special case of unweighted $\ell^1$ minimization, our guarantees reduce to those introduced previously in \cite{HamptonDoostanCSPCE,Rauhut,YanGuoXui_l1UQ} (although our theorems avoid the issues surrounding  infinite expansions mentioned above).  For weighted $\ell^1$ minimization, a corollary of our main result yields guarantees similar to those in \cite{RauhutWardWeighted}.  As we demonstrate, however, the guarantees of \cite{RauhutWardWeighted} are not sharp for a large class of functions and weights.  Fortunately, by returning to our abstract result we derive an improved recovery guarantee which is sharp for this class.  

We next use our main result to identify a good overall weighting strategy, in the sense that minimizes the number of measurements required in our recovery guarantee.  Similarly to recent work \cite{ChkifaDownwardsCS}, we show that this choice of weights yields a reconstruction procedure for tensor Chebyshev and Legendre polynomial expansions that mitigates the curse of dimensionality to a substantial extent.  This is on the proviso that the coefficients of the function being approximated possess a particular type of structured sparsity defined by so-called \textit{lower sets}, which is a reasonable assumption in applications of interest (see \cite{ChkifaEtAl,ChkifaDownwardsCS,CohenDeVoreSchwabFoCM,CohenDeVoreSchwabRegularity,MiglioratiJAT,MiglioratiEtAlFoCM} and references therein).  Finally, we apply our abstract recovery guarantee to assess an alternative strategy where the weights are chosen based on prior support information; a strategy which has been advocated in a number of recent works \cite{PengHamptonDoostantweighted,KarniadakisUQCS}.  When the weights our chosen in this way, our recovery guarantee provides a theoretical basis for the empirically-observed benefits of this approach.

The outline of the remainder of this paper is as follows.  In \S \ref{s:background} we present relevant background material on CS for function approximation, and give an overview of our main contributions in this paper.  We give some preliminary notation in \S \ref{s:preliminaries}.  In \S \ref{s:infdiml1} we introduce the infinite-dimensional weighted $\ell^1$ minimization problem, and in \S \ref{ss:main_examp} we present the main examples that will be used to demonstrate the our results.  Our main abstract recovery guarantee is presented in \S \ref{ss:main_res}, and in \S \ref{s:consequences} we discuss various consequences of it, including its application to high-dimensional approximation using polynomials.  Finally, in \S \ref{s:proof} we present the proof of our main result.

\section{Background and overview}\label{s:background}
In this section, we first give an overview previous work on CS for function approximation and then present a summary of the main results of this paper.

\subsection{Previous work}
In recent years, numerous works have sought to apply the techniques of CS to multivariate function approximation.  
Theoretical guarantees for recovery of sparse polynomial expansions via unweighted $\ell^1$ minimization were first presented in \cite{Rauhut} for the univariate case and \cite{YanGuoXui_l1UQ} for the multivariate case.  Weighted $\ell^1$ minimization was proposed in \cite{PengHamptonDoostantweighted,KarniadakisUQCS}, and theoretical results based on weighted sparsity were presented in \cite{RauhutWardWeighted}.  For theoretical results on lower sets related to this paper, see \cite{ChkifaDownwardsCS}.  Applications of CS techniques to uncertainty quantification have been pursued numerous works, including \cite{DoostanOwhadiSparse,MathelinGallivanCSPDErandom,PengHamptonDoostantweighted,KarniadakisUQCS} and references therein.  Besides random sampling (from appropriate continuous measures) several works have also proposed new sampling strategies for CS that aim to improve reconstruction quality.  These include coherence-based sampling \cite{HamptonDoostanCSPCE}, preconditioning \cite{JakemanEtAlChristoffel,NarayanJakemanZhouChristoffelLS}, deterministic sampling \cite{XuZhouSparseDeterministic} and subsampling from deterministic Gaussian quadratures \cite{TangIaccarino}.    Outside of CS theory, worst-case recovery guarantees for weighted $\ell^1$ minimization for deterministic sampling were presented in \cite{Adcockl1pointwise}, demonstrating near-optimal performance for general scattered data.

\subsection{Compressive function approximation}\label{ss:current_approaches}
Let $\{ \phi_i \}_{i \in \bbN}$ an orthonormal basis of functions (e.g.\ polynomials) and consider a function $f = \sum_{i \in \bbN} x_i \phi_i$.  Suppose that $\{ t_i \}^{m}_{i=1}$ is a finite set of points, typically chosen randomly from an appropriate distribution, and consider the measurements $y = \{ f(t_i) \}^{m}_{i=1}$.  In order to approximate $f$, it suffices to approximate its coefficients $x$ from the measurements $y$.  However, $x$ is typically an infinite vector, meaning that some sort of discretization is required.  In a majority of previous works, this discretization is performed first; that is, prior to formulating the optimization problem.  Specifically, 
one introduces a fixed $N \geq m$ and seeks to approximate the first $N$ coefficients $x_1,\ldots,x_N$ of $x$ by solving the following inequality-constrained (weighted) $\ell^1$ minimization problem:
\be{
\label{intro_bad_min}
\min_{z \in \bbC^N} \| z \|_{1,w}\ \mbox{subject to $\| A z - y \| \leq \delta$}.
}
Here $\delta \geq 0$, $\| z \|_{1,w} = \sum^{N}_{i=1} w_i | z_i|$ is the $\ell^1_w$-norm on $\bbC^N$ with weights $w_i >0$ (we discuss the issue of weights next) and $A = \{ \phi_j(t_i) \}^{m,N}_{i=1,j=1} \in \bbC^{m \times N}$.  The parameter $\delta$ is an artefact of the discretization, and is chosen so that the exact coefficients $\{ x_i \}^{N}_{i=1}$ are feasible for \R{intro_bad_min}, i.e.\
\be{
\label{finite_bad}
\nm{ f - \sum^{N}_{i=1} x_i \phi_i }_{L^\infty} \leq \delta.
}
If $\hat{x}$ is a minimizer of \R{intro_bad_min}, then one defines the corresponding approximation to $f$ as $\tilde{f} = \sum^{N}_{i=1} \hat{x}_i \phi_i$.

As discussed in \cite{Adcockl1pointwise}, the error committed by the approximation $\tilde{f}$ will depend on the choice of $\delta$.  Hence a good estimation of the norm of the expansion tail is important to ensure accurate results \cite{KarniadakisUQCS}.  Herein lies a problem.  In general, this tail error is unknown.  Whilst techniques such as cross validation \cite{DoostanOwhadiSparse,PengHamptonDoostantweighted,KarniadakisUQCS} have been used to provide practical estimations for $\delta$, these are both computationally expensive and wasteful in terms of the data.  A key element of the framework we develop in this paper is that it does not require knowledge of $\delta$, and thus allows one to avoid this estimation step.

In addition to this practical issue, from a theoretical perspective all existing CS recovery guarantees for \R{intro_bad_min} (see \cite{ChkifaDownwardsCS,DoostanOwhadiSparse,HamptonDoostanCSPCE,RauhutWardSpher,Rauhut,RauhutWardWeighted,YanGuoXui_l1UQ}) assume \textit{a priori} knowledge of $\| f - \sum^{N}_{i=1} x_i \phi_i \|_{L^\infty}$.  Whilst this is not a problem when $f$ is itself a polynomial of known degree $N$, such guarantees are less informative about the approximation capabilities of weighted $\ell^1$ minimization for objects with infinite expansions, i.e.\ functions.  The theoretical guarantees we present in this paper seek to overcome this limitation.

\subsection{Weighted $\ell^1$ minimization}\label{ss:weights}
In a number of recent works it has been observed empirically that unweighted $\ell^1$ minimization (i.e.\ with $w_i = 1$, $\forall i$) often gives relatively poor approximations, and that better results are possible if slowly growing weights $w_i$ are introduced \cite{PengHamptonDoostantweighted,RauhutWardWeighted,KarniadakisUQCS}.  For deterministic samples $\{ t_i \}^{m}_{i=1}$, this was explained recently in \cite{Adcockl1pointwise}.  Therein it was shown that unweighted $\ell^1$ minimization is in general unsuitable for function approximation, since it suffers from a so-called \textit{aliasing} phenomenon stemming from the infinite-dimensionality of the problem.  In the unweighted case, there are solutions $\hat{x} = \{ \hat{x}_i \}^{\infty}_{i=1}$ to the optimization problem with nonzero coefficients far out in the expansion tail.  The corresponding approximation $\tilde{f} = \sum_{i \in \bbN} \hat{x}_i \phi_i$ fits the data, but oscillates rapidly in between the data points, due to nonzero high-frequency modes.  The introduction of slowly growing weights removes this phenomenon, however, since high-frequency modes are penalized by increasing weights.

\begin{figure}[t]
\begin{center}
\begin{tabular}{ccc}
\includegraphics[width=5.25cm]{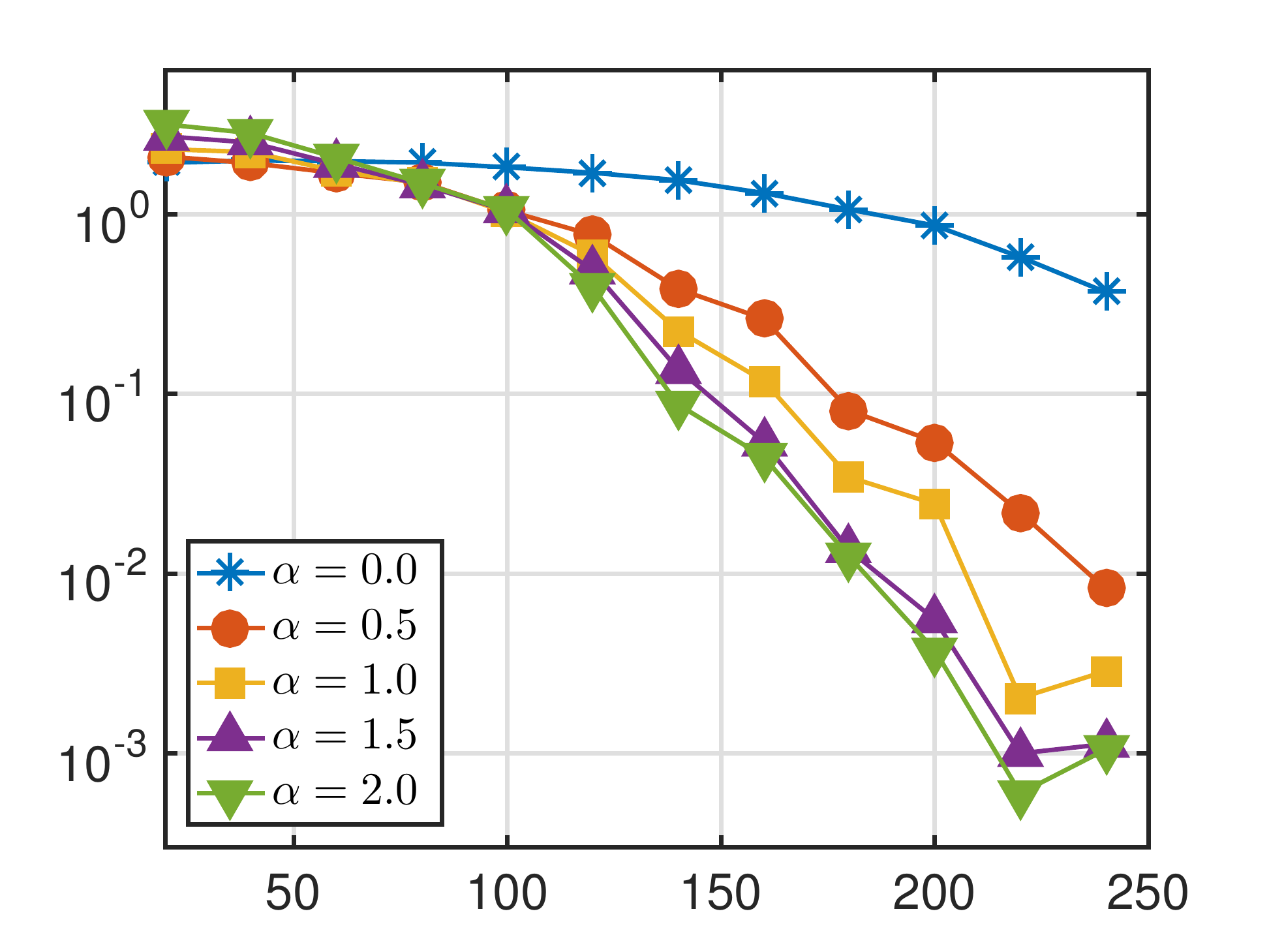} &
\includegraphics[width=5.25cm]{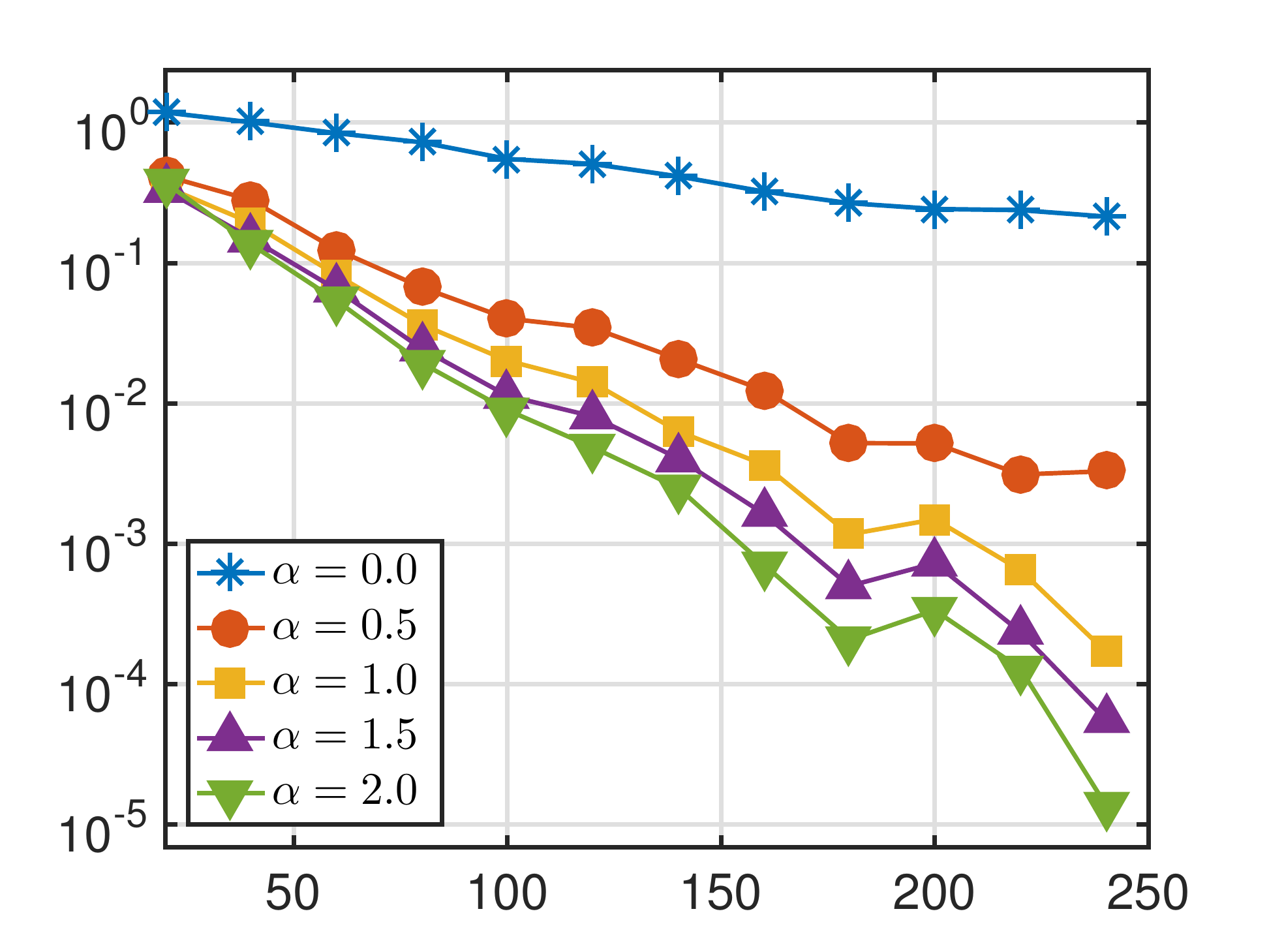} &
\includegraphics[width=5.25cm]{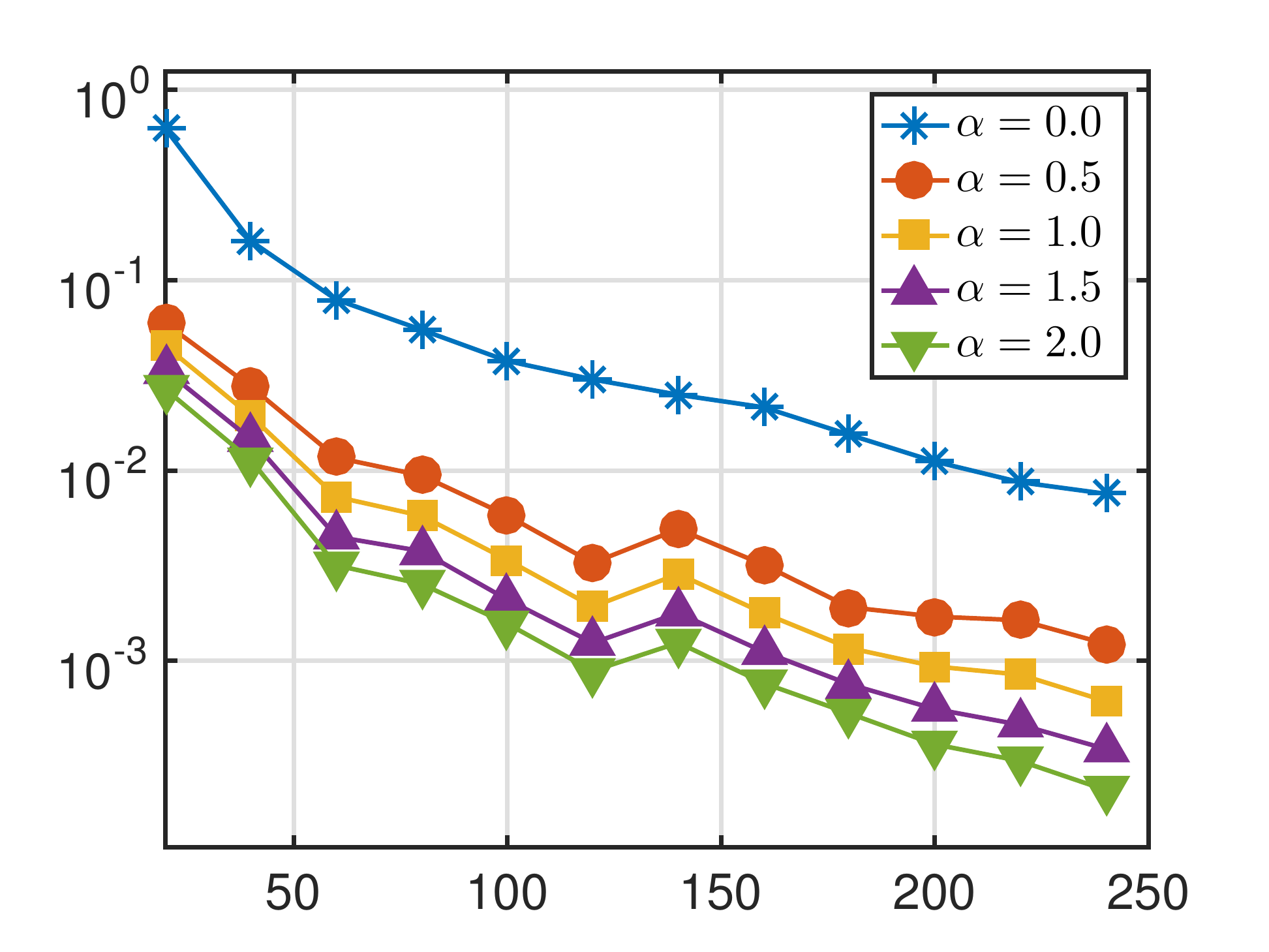} 
\\
\includegraphics[width=5.25cm]{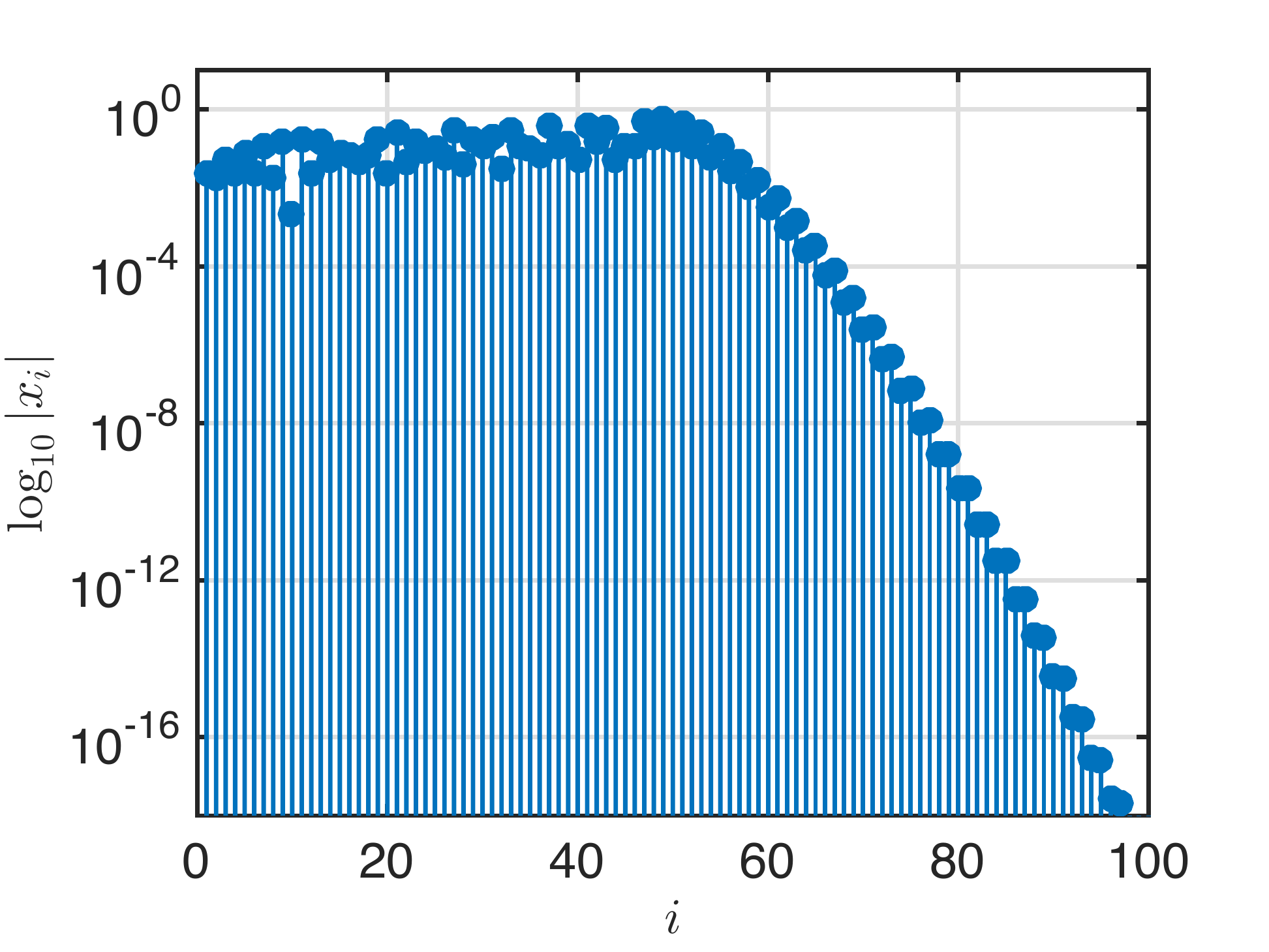} &
\includegraphics[width=5.25cm]{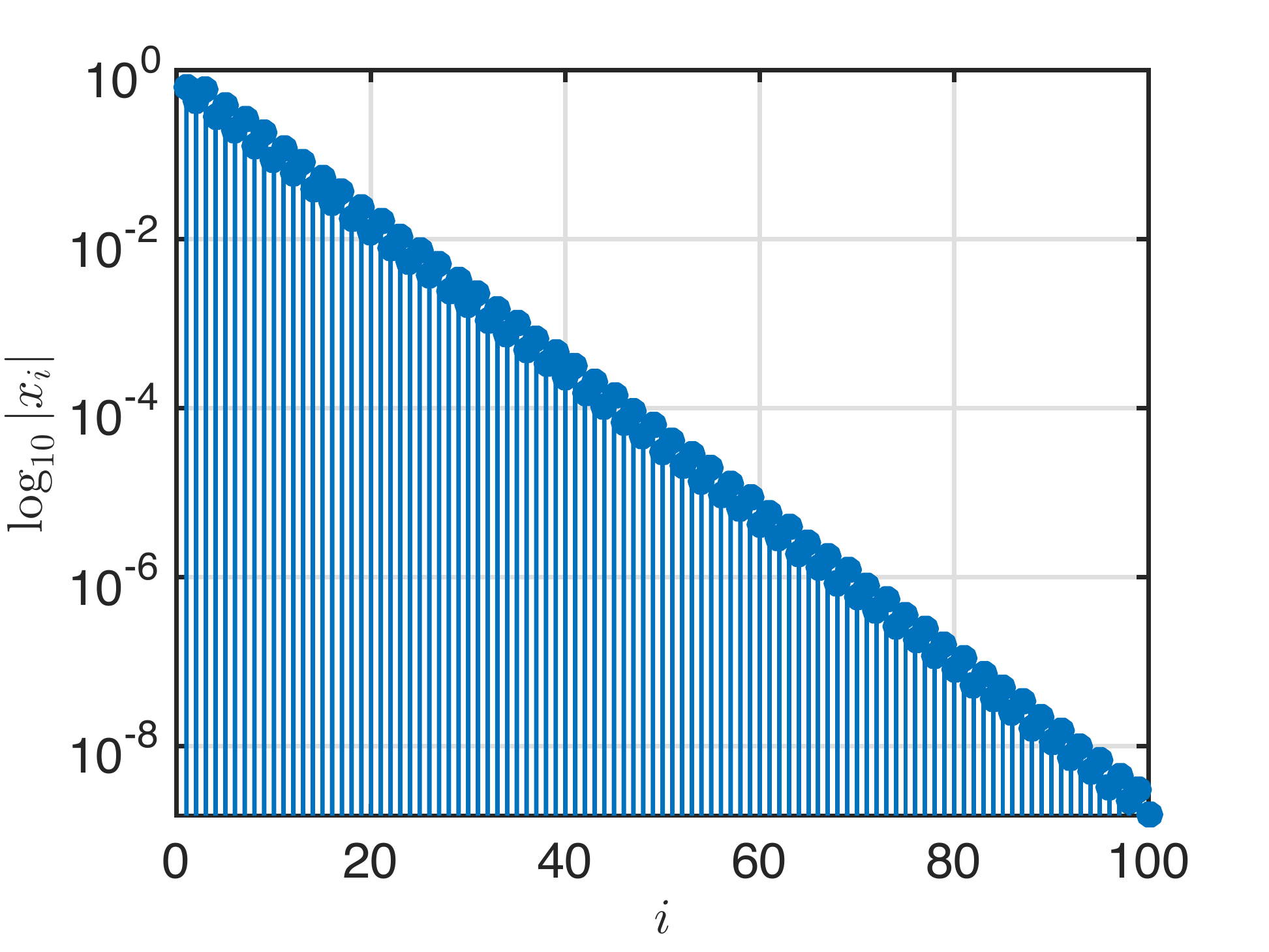} &
\includegraphics[width=5.25cm]{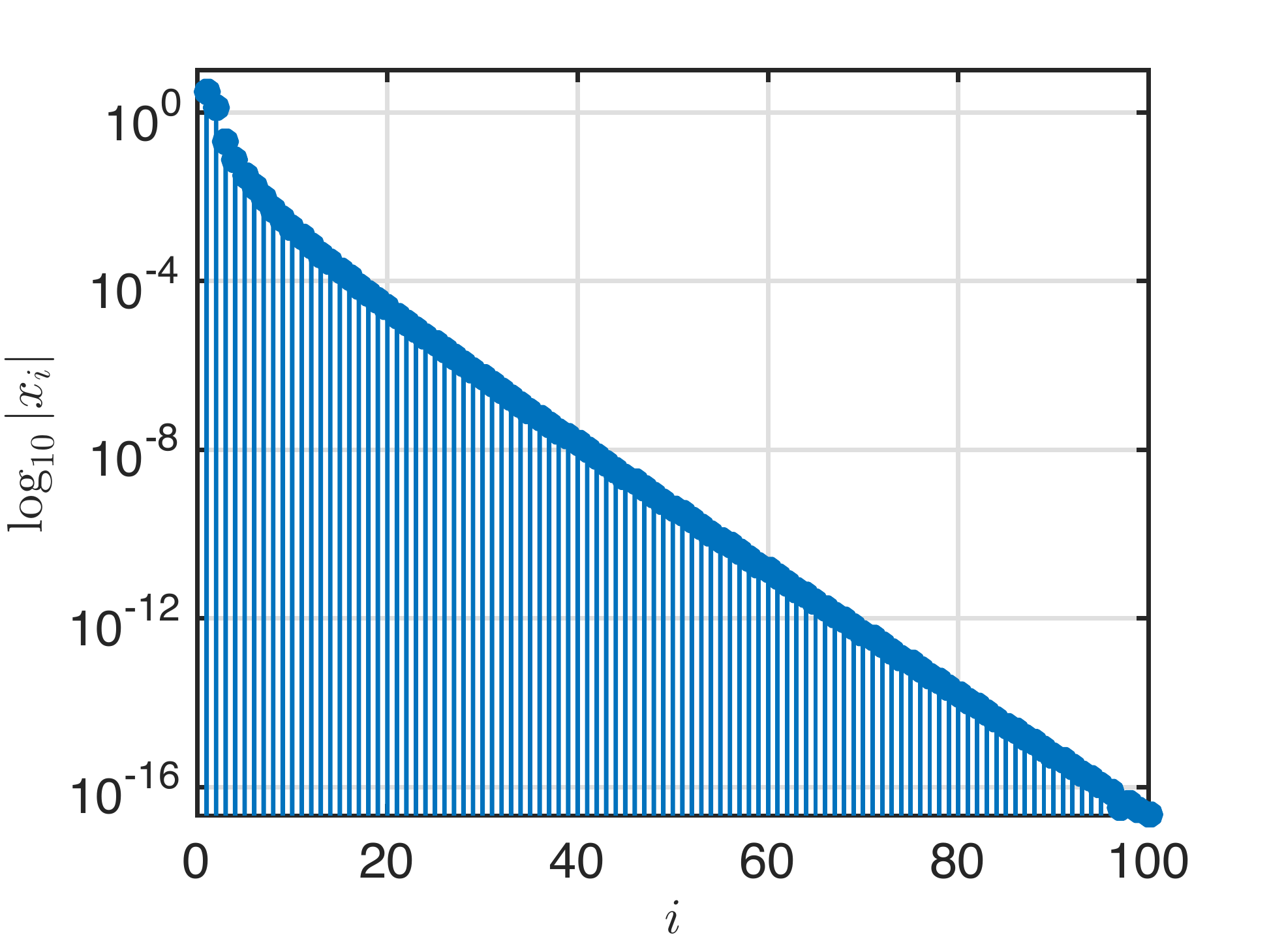} 
\\
CC, $f(t) = \cos(36 \sqrt{2} t + 1/3)$  & LC, $f(t) = \frac{1+3t}{1+50 t^2}$ & LU, $f(t) = \sqrt{1.05+t}$
\end{tabular}
\end{center}
\caption{
Top row: The error $\| f - \tilde{f} \|_{L^\infty}$ (averaged over $50$ trials) against $m$ for Chebyshev (C) or Legendre (L) polynomials with points drawn from the Chebyshev (C) or uniform (U) measure.  Here $\tilde{f} = \sum_{i \in I_K} \hat{x}_i \phi_i$, where $\hat{x}$ is a solution of \R{fin_min} and $I_K = \{0,\ldots,K\}$.  The parameter $K = 1000$ was used and the weights were taken to be $w_i = (i+1)^{\alpha}$ for various $\alpha \geq 0$.  Bottom row: Chebyshev or Legendre coefficients of the function $f(t)$.  In this and all other experiments in this paper, we use the SPGL1 package \cite{spgl1:2007,BergFriedlander:2008} with a maximum of 100,000 iterations.
}
\label{f:ErrIntroduction}
\end{figure}

The focus of this paper is on samples $\{ t_i \}^{m}_{i=1}$ drawn randomly from an appropriate distribution.  Unlike for deterministic samples, approximations computed by solving the unweighted $\ell^1$ minimization problem do converge (with high probability) as the number of samples increases.  Yet adding weights generally leads to a more rapidly-decreasing approximation error in this setting \cite{PengHamptonDoostantweighted,KarniadakisUQCS,RauhutWardWeighted}; see Fig.\ \ref{f:ErrIntroduction} for an illustration.  A central contribution of this paper is to understand and quantify this benefit theoretically.  Our basis for this is a new recovery guarantee for weighted $\ell^1$ minimization.

\subsection{Recovery guarantees}\label{ss:recov_intro}
Standard CS theory states that one can recover a vector $x$ of sparsity $s$, i.e.\
\be{
\label{sparsity}
s= | \Delta | = \sum_{i \in \Delta } 1,\qquad \Delta = \{ i : x_i \neq 0 \},
}
using, up to log factors, $m \approx s$ appropriately-chosen measurements, regardless of the locations of the nonzero entries of $x$.  In practice, this can be achieved by solving an $\ell^1$ minimization problem.  When considering weighted $\ell^1$ minimization on the other hand, it was proposed in \cite{RauhutWardWeighted} to replace sparsity \R{sparsity} by the following weighted sparsity measure
\bes{
s = |\Delta |_w =  \sum_{i \in \Delta } w^2_i,\qquad \Delta = \{ i : x_i \neq 0 \},
}
corresponding to the optimization weights $w = \{w_i \}_{i \in \bbN}$.  The work of \cite{RauhutWardWeighted} has established a measurement condition of the form
\be{
\label{current_weighted_main} 
m \approx | \Delta |_{w} \times \mbox{log factors},
}
for weighted $\ell^1$ minimization with appropriate measurements.

Unfortunately, such guarantees do not explain the observed empirical performance of weighted $\ell^1$ minimization in some important cases.  For example, suppose that $\Delta = \{1,\ldots,M\}$.  That is, $x$ is nonzero in its first $M$ entries, or more generally (for inexact sparsity), the largest $M$ entries of $x$ are also its first $M$ entries.  This phenomenon is common in polynomial expansions of smooth functions, especially in lower-dimensional settings.  Herein coefficients tend to exhibit decay (see Fig.\  \ref{f:ErrIntroduction}), with the largest $M$ coefficients often coinciding with the first $M$, or more generally, the first $\ord{M}$.  A simple example of this is an oscillatory function, for which the coefficients $x_i$ are $\ord{1}$ up to a certain resolution criterion, after which they are numerically zero (see Fig.\  \ref{f:ErrIntroduction}).

Suppose for simplicity that the weights $w_i = i^{\alpha}$ are taken to be polynomially growing with index $\alpha$.  Then the recovery guarantee \R{current_weighted_main} reduces to
\be{
\label{current_weighted}
m \approx M^{2\alpha+1} \times \mbox{log factors}.
}
In other words, more rapidly-growing weights seemingly require more measurements to recover $x$.  However, this conclusion is at odds with numerical results shown in Fig.\  \ref{f:ErrIntroduction}, wherein it is observed that increasing the weights in fact decreases the error to a moderate degree, and certainly does not worsen it.  This example is representative of quite general behaviour of weighted $\ell^1$ minimization.

Fortunately, it transpires that \R{current_weighted} is not sharp.  A corollary of our main result gives an estimate for this example which is both independent of the weights, in certain cases, and provably sharp.

\subsection{Contributions}\label{ss:contributions}
Our main contribution is a general result on infinite-dimensional CS for function interpolation based on weighted $\ell^1$ minimization.  A simplified version (see Remark \ref{r:simple}) of our main result is as follows:

\thm{
\label{t:simple_main}
Let $0 < \epsilon < \E^{-1}$, $D \subseteq \bbR^d$ be a domain with a probability measure $\nu$, $\{ \phi_i \}_{i \in \bbN}$ be an orthonormal system in $L^2_\nu(D)$ with $u_i =  \| \phi_i \|_{L^\infty} < \infty$ for all $i \in \bbN$, $w = \{ w_i \}_{i \in \bbN}$, $w_i > 0$ be weights and $f = \sum_{i \in \bbN} x_i \phi_i$ with $x = \{ x_i \}_{i \in \bbN} \in \ell^1_{w}(\bbN)$.  Let $K \in \bbN$ and $\Delta \subseteq \{1,\ldots,K\}$ be a subset with $\min_{i \in \{1,\ldots,K\} \backslash \Delta } \{ w_i \} \geq 1$ and suppose that $t_1,\ldots,t_m$ are drawn independently from the measure $\nu$, where 
\be{
\label{simple_guarantee}
m \gtrsim \left ( | \Delta |_u + \sup_{i \in \{1,\ldots,K\} \backslash \Delta} \{ u^2_i/w^2_i \} | \Delta |_w \right ) \cdot \log(\epsilon^{-1}) \cdot \log(| \Delta |_w K).
}
Then, by solving a weighted $\ell^1$ minimization problem of size $m \times K$ with weights $w$, it is possible to approximate $x$ (and therefore $f$) from the data $y = \{ f(t_i) \}^{m}_{i=1}$ up to an error proportional to
\bes{
\sum_{i \notin \Delta} w_i |x_i| + T_{K,w}(x),
}
with probability $1-\epsilon$, where $T_{K,w}(x)$ is the truncation error \R{TK_def}.  Furthermore, in the absence of noise the resulting approximation interpolates the data $\{f(t_i)\}^{m}_{i=1}$.
}

This result highlights the main contributions of this paper:

\vspace{1pc} \noindent 1.\ \textit{Removal of a priori tail estimates.}\  Theorem \ref{t:simple_main}, unlike recovery guarantees based on \R{intro_bad_min}, completely avoids the need for \textit{a priori} knowledge of the expansion tail.  In particular, it applies to arbitrary functions, not just finite expansions in the system $\{ \phi_i \}_{i \in \bbN}$.  
This result of this is the additional truncation term $T_{K,w}(x)$ in the error estimate.  But, as we explain in \S \ref{ss:main_res}, this term can be estimated and is typically negligible for large enough $K$.  

\vspace{1pc} \noindent 2.\ \textit{A weights-independent recovery guarantee.}\  As in \S \ref{ss:recov_intro}, let $\Delta = \{1,\ldots,M\}$ and $w_i = i^{\alpha}$.  It is reasonable to assume that $u_i = \ord{i^\beta}$ as $i \rightarrow \infty$ for some $\beta > 0$.  If $\alpha \geq \beta$, then \R{simple_guarantee} reduces to 
\bes{
m \gtrsim M^{2\beta+1} \cdot \log(\epsilon^{-1}) \cdot \log(M^{2\alpha+1} K).
}
In contrast to \R{current_weighted}, this condition is independent of the weights parameter $\alpha$.  In particular, for the recovery of Chebyshev polynomial expansions, where $\beta = 0$, it reduces to $m = \ord{M}$ (up to log factors) and for Legendre polynomial expansions, in which case $\beta = 1/2$, we obtain $m=\ord{M^2}$.  Both results are essentially sharp:

\rem{
\label{r:quadscalingsharp}
Clearly the result for Chebyshev polynomials is sharp up to log factors, since it is linear in the number of unknowns $M$.  In \cite{Adcockl1pointwise}, using results from \cite{AdcockNecSamp,TrefPlatteIllCond}, it was shown that no robust method can recover the first $M$ Legendre polynomial coefficients using asymptotically fewer than $m \asymp M^2$ measurements when the samples are exactly equidistributed (as opposed to randomly drawn) according to the uniform measure.
}

\vspace{1pc}\noindent 3.\ \textit{Identifying good weights.}\ Log factors aside, the number of measurements stipulated by \R{simple_guarantee} is dependent on the factor
\bes{
 | \Delta |_u + \sup_{i \in \{1,\ldots,K\} \backslash \Delta} \{ u^2_i/w^2_i \} | \Delta |_w.
}
Note that the first term is independent of the optimization weights, and depends only on the the \textit{intrinsic} weights $u$.  Therefore, in the absence of any further assumptions on $\Delta$ (see below for this case), a good weighting strategy would seek to make the second term equal to the first.  This can be achieved by setting
\be{
\label{gen_opt_weights}
w_i = u_i,\quad \forall i,
}
which leads to the following recovery guarantee for $\ell^1_u$ minimization:
\bes{
m \gtrsim | \Delta |_u \cdot \log(\epsilon^{-1}) \cdot \log(|\Delta|_u K).
}
Note that, beside constants and log factors, this is strictly smaller than the recovery guarantee 
\bes{
m \gtrsim \left ( | \Delta |_u + \max_{i=1,\ldots,K} \{ u^2_i \} | \Delta | \right ) \cdot \log(\epsilon^{-1}) \cdot \log(|\Delta|_u K),
}
for unweighted $\ell^1$ minimization, and may be substantially smaller depending on the support of $\Delta$.

\vspace{1pc} \noindent 4.\ \textit{Overcoming the curse of dimensionality.}\  In the case of polynomial approximations in high dimensions, the weighting strategy \R{gen_opt_weights} leads to a substantially smaller recovery guarantee for certain structured sparse support sets $\Delta$; specifically, so-called lower sets (see Definition \ref{d:lower_set}).  These sets are known to be good models for the support sets of polynomial coefficients in high dimensions (see \cite{ChkifaEtAl,ChkifaDownwardsCS,CohenDeVoreSchwabFoCM,CohenDeVoreSchwabRegularity,MiglioratiJAT,MiglioratiEtAlFoCM} and references therein).  For example, with tensor Chebyshev polynomials and points drawn randomly from the tensor Chebyshev measure the recovery guarantee \R{simple_guarantee} yields
\bes{
m \gtrsim 2^d s \times \mbox{log factors},\qquad \mbox{if $w_i = 1$, $\forall i$},
} 
where $s = |\Delta|$.  However, if $\Delta$ is also a lower set, one has
\bes{
m \gtrsim s^{\log(3)/\log(2)}  \times \mbox{log factors},\qquad \mbox{if $w_i = u_i$, $\forall i$}.
}
Log factors aside (see Remarks \ref{r:logfactors} and \ref{r:logfactorslower}), the latter does not suffer from the curse of dimensionality and agrees with the best known estimates for least-squares approximation in lower sets \cite{ChkifaEtAl}.  Similarly, for tensor Legendre polynomials with points drawn randomly from the uniform measure, the recovery guarantee \R{simple_guarantee} yields
\bes{
m \gtrsim s^2 \times \mbox{log factors},\qquad \mbox{if $w_i = u_i$, $\forall i$},
}
for lower sets $\Delta$ (this result is essentially sharp -- recall Remark \ref{r:quadscalingsharp}), whereas the unweighted guarantee grows exponentially with the dimension of the truncated polynomial space.  We refer to \S \ref{ss:recovery_examples} for further details.

\vspace{1pc} \noindent 5.\ \textit{Support estimation via weighted $\ell^1$ minimization.}\  Since our recovery guarantee applies to any choices of weights $w_i$, it can be used to explain why the strategy of choosing weights based on \textit{a priori} knowledge of part of the support set of the coefficients reduces the number of measurements required.  Our main result in this setting, detailed in \S \ref{ss:support_est}, shows that if over half of the support set is correctly estimated, then such a weighting strategy leads to a strictly smaller recovery guarantee than that of the unweighted case.

\rem{
\label{r:simple}
Theorem \ref{t:simple_main} makes several simplifications for illustrative purposes that are not necessary in our main result, Theorem \ref{t:main}.  First, the measurements $t_1,\ldots,t_m$ can be drawn from a different measure $\mu$ than the orthogonality measure $\nu$ of the functions $\{ \phi_i \}_{i \in \bbN}$.  Second, the measurements can be noisy.  In our main result we also allow for indexing over arbitrary countable sets $I$, which is particularly useful in the case of multivariate polynomial approximations.
}

\section{Preliminaries}\label{s:preliminaries}
Throughout this paper, $D \subseteq \bbR^d$ will be a domain and $\nu$ will be an integrable nonnegative weight function on $D$.  We write $L^2_{\nu}(D)$ for the space of complex-valued weighted square-integrable functions on $D$, with norm $\nm{\cdot}_{L^2_{\nu}}$ and inner product $\ip{\cdot}{\cdot}_{L^2_{\nu}}$.  Given $D$ and $\nu$, we let $\{ \phi_i \}_{i \in I} \subseteq L^2_{\nu}(D) \cap L^\infty(D)$ be a set of functions that are orthonormal with respect $\nu$, where $I$ is an index set that is at most countable.  For a function $f \in \overline{\spn\{\phi_i : i \in I\}}$ we let $x = \{ x_i \}_{i \in I} \in \ell^2(I)$ be its coefficients in the system $\{ \phi_i \}_{i \in I}$, i.e.\
\bes{
f = \sum_{i \in I} x_i \phi_i,\qquad x_i = \ip{f}{\phi_i}_{L^2_{\nu}}.
}
Our goal throughout this paper is to recover the coefficients $x$ from a small number of pointwise evaluations of $f$.

Several other pieces of standard notation will also be used.  The space of square-summable sequences indexed over $I$ will be denoted by $\ell^2(I)$, with its norm and inner product given by $\nm{\cdot}$ and $\ip{\cdot}{\cdot}$ respectively.  The set $\{ e_i \}_{i \in I}$ denotes the canonical basis of $\ell^2(I)$ and, if $\Delta \subseteq I$, we write $P_{\Delta}$ for the orthogonal projection onto $\overline{\spn\{e_i : i \in \Delta \}}$.  Whenever convenient, we will also regard $P_\Delta$ for finite $\Delta$ as a mapping with range $\bbC^{|\Delta|}$.  Similarly, we will consider a vector $x \in \mathrm{Ran}(P_\Delta)$ interchangeably as a sequence $x \in \ell^2(I)$ supported on the set $\Delta$ and as a vector in $\bbC^{|\Delta|}$.

\subsection{Weighted spaces and sparsity}
For the remainder of this paper, $w = \{ w_i \}_{i \in I}$ will be a set of positive weights used in the optimization problem.  Define the space of weighted summable sequences by
\bes{
\ell^1_w(I) = \left \{ x = \{ x_i \}_{i \in I} : \| x \|_{1,w} : = \sum_{i \in I} w_i | x_i | < \infty \right \}.
}
Note that we may write the norm as $\| x \|_{1,w} = \| W x \|_{1}$, where 
\be{
\label{W_def}
W = \mathrm{diag}(w_1,w_2,\ldots),
}
is the infinite diagonal matrix of weights and $\nm{\cdot}_{1}$ is the standard unweighted $\ell^1$ norm.  For a set $\Delta \subseteq I$ we also define its weighted cardinality by
\bes{
| \Delta |_{w} = \sum_{i \in \Delta} (w_i)^2.
}
If $w_i = 1$, $\forall i \in I$, then we merely write $| \Delta |$ for the corresponding unweighted cardinality of $\Delta$.

\subsection{Sampling points and the operator $U$}
Given $\nu$ and $\{ \phi_i \}_{i \in I}$ let $\mu$ be a probability measure on $D$ satisfying
\be{
\label{mu_cond}
\sup_{t \in D} \sqrt{\nu(t)/\mu(t)} | \phi_i(t) | < \infty,\quad \forall i \in I,
}
and define
\be{
\label{u_i_def}
u_i : =  \max \left \{ 1, \sup_{t \in D} \sqrt{\nu(t)/\mu(t)} | \phi_i(t) |  \right \},\quad \forall i \in I.
}
Note that we do not require the sequence $u = \{u_i \}_{i \in I}$ to be uniformly bounded in $i$.  Given such a probability measure $\mu$, we assume from now on that the sampling points $t_1,\ldots,t_m$ are drawn independently from $\mu$.

With $\nu$, $\{ \phi_i \}_{i \in I}$, $\mu$ and $\{ t_i \}^{m}_{i=1}$ in hand, we define the infinite matrix $U$ as follows:
\be{
\label{U_def}
U = \left \{ \phi_j(t_i) \sqrt{\nu(t_i) / \mu(t_i) } \right \}^{m,\infty}_{i=1,j=1}.
}
We shall view $U$ interchangeably as an infinite matrix and also an operator.  The following lemma identifies an appropriate domain for $U$ to be bounded:

\lem{
\label{l:U_bounded}
Let $u = \{u_i \}_{i \in I}$ and $U$ be as in \R{u_i_def} and \R{U_def} respectively.  Then the operator $U : \ell^1_u(I) \rightarrow \bbC^m$ is bounded.
}
\prf{
Let $x \in \ell^1_u(I)$.  Then $| (U x)_i | = \sum_{j \in I} | x_j | |\phi_j(t_i) | \sqrt{\nu(t_i) / \mu(t_i) } \leq \sum_{j \in I} u_j |x_j | = \| x \|_{1,u}$.
}

\section{Infinite-dimensional weighted $\ell^1$ minimization}\label{s:infdiml1}

In this section, we introduce the infinite-dimensional weighted $\ell^1$ minimization formulation that will be used throughout this paper.  

Given a function $f$ and sampling points $\{ t_i \}^{m}_{i=1}$, the measurements will be of the form
\bes{
 f(t_i ) + \tilde{e}_i,\quad i=1,\ldots,m,
}
where $\tilde{e}_i$ are noise terms satisfying the weighted estimate
\bes{
\sum^{m}_{i=1} \frac{\nu(t_i)}{\mu(t_i)} | \tilde{e}_i |^2 \leq \eta^2,
}
for some known noise parameter $\eta$.  Define the scaled noise vector
\bes{
e = \{ e_i \}^{m}_{i=1},\quad e_i = \sqrt{\nu(t_i)/\mu(t_i)} \tilde{e_i},\qquad \| e \| \leq \eta,
}
and note that the measurements $y = \{y_i \}^{m}_{i=1}$ can be expressed as
\bes{
y = U x + e,\qquad y_i = \left \{ \sqrt{\nu(t_i)/\mu(t_i)} f(t_i) + e_i \right \}^{m}_{i=1},
}
where $U$ is as in \R{U_def}.  Suppose now that $w = \{w_i \}_{i \in I}$ are weights.  As in \cite{Adcockl1pointwise}, we consider the optimization problem
\be{
\label{inf_min}
\inf_{z \in \ell^1_w(I)} \| z \|_{1,w}\ \mbox{subject to $\| U z - y \| \leq \eta$.}
}
Note that in the absence of noise, i.e.\ $\eta = 0$, solutions of the problem exactly interpolate the function $f$ at the points $t_1,\ldots,t_m$.  Moreover, unlike \R{intro_bad_min} there is no need to know bounds for the expansion tail in order to formulate \R{inf_min}.

Unfortunately, this problem cannot be solved numerically, since it involves minimizing over an infinite-dimensional space.  To overcome this, we need to truncate \R{inf_min} in such a way so that we retain the important properties of \R{inf_min} noted above.  We do this as follows.  For $K=1,2,\ldots$ let $I_K \subseteq I$ be a subset of $I$ of finite cardinality and let $P_K = P_{I_K}$ denote the projection onto $I_K$.  We shall assume that the sequence of projections $\{P_K\}_{K \in \bbN}$ converges strongly to the identity operator on $\ell^2(I)$, i.e.\
\bes{
P_{K} x \rightarrow x,\quad \forall x \in \ell^2(I).
}
Note that we do not require the sets $I_K$ to be nested, although this will often be the case in practice.  Given the subsets $I_K$, we now replace \R{inf_min} with the following problem:
\be{
\label{fin_min}
\min_{z \in P_K(\ell^1_w(I))} \| z \|_{1,w}\  \mbox{subject to $\| U P_K z - y \| \leq \eta$}.
}
This problem is equivalent to a minimization problem on $\bbC^K$, and therefore numerically solvable.  In particular, the space $P_K(\ell^1_w(I))$ is isomorphic to $\bbC^K$ and $U P_K$ is equivalent to an $m \times K$ matrix formed by the columns of $U$ with indices in $I_K$. Importantly, however, \R{fin_min} retains the key features of \R{inf_min}.  Namely, there is no need to know the expansion tail, and in the absence of noise solutions of \R{fin_min} interpolate $f$ exactly at the data points.

A key issue for \R{fin_min} is the choice of the \textit{truncation parameter} $K$.  Loosely speaking, $K$ should be taken sufficiently large such that the additional error induced by solving \R{fin_min} instead of \R{inf_min} is small.  It is important, however, that this criterion be independent of the function $f$ to approximate, i.e.\ it should not involve the expansion tail $\sum_{i \in I \backslash I_K} x_i \phi_i$, since \textit{a priori} estimates for this term are generally unknown (recall \S \ref{ss:current_approaches}).  Fortunately, as we demonstrate in \S \ref{ss:main_res}, this is indeed possible.

\section{Main examples: tensor Chebyshev and Legendre polynomials}\label{ss:main_examp}
To illustrate the main results of this paper, we consider the case of tensor products of Chebyshev and Legendre polynomials on the unit hypercube $D = (-1,1)^d$.  Recall that one-dimensional Legendre and Chebyshev polynomials are orthogonal with respect to the measures
\bes{
\nu(t) = 1/2 \quad \mbox{and}\quad \nu(t) = \frac{1}{\pi \sqrt{1-t^2}},\qquad t \in (-1,1),
}
respectively.  In the hypercube, the corresponding orthogonality measures are
\bes{
\nu(t) = 2^{-d}\quad \mbox{and} \quad \nu(t) = \prod^{d}_{j=1} \frac{1}{\pi(1-t^2_j)^{1/2}},\qquad t = (t_1,\ldots,t_d) \in (-1,1)^d.
}
If $\phi_0,\phi_1,\ldots$ are the univariate polynomials of Chebyshev or Legendre type, we define the multivariate system via tensor products:
\bes{
\phi_i(t) = \prod^{d}_{j=1} \phi_{i_j}(t_j),\qquad t = (t_1,\ldots,t_d) \in (-1,1)^d,\  i = (i_1,\ldots,i_d) \in \bbN^d_0.
}
Within this setup, we shall address the following three specific sampling scenarios, all of which are permissible under the condition \R{mu_cond}:

\vspace{1pc} \noindent \textit{Tensor Chebyshev polynomials, random sampling from the Chebyshev measure.}  In this case, the measures are given by $\nu(t) = \mu(t) = \prod^{d}_{j=1} \frac{1}{\pi (1-t^2_j)^{1/2}}$.  Since univariate Chebyshev polynomials satisfy $| \phi_0(t) | =1$ and $\sup_{t \in (-1,1)} |\phi_i(t) | = \sqrt{2}$ otherwise we find that
\be{
\label{CC_ui}
u_i = \sup_{t \in D} | \phi_i(t) | = 2^{|i|_0/2},
}
for this example, where $|i|_0 = | \{ j : i_j \neq 0 \} |$ for $i = (i_1,\ldots,i_d) \in \bbN^d_0$.

\vspace{1pc}\noindent  \textit{Tensor Legendre polynomials, random sampling from the uniform measure.}  In this case, $\nu(t) = \mu(t) = 2^{-d}$.  Since the univariate Legendre polynomial satisfies $\sup_{t \in (-1,1)} | \phi_i(t) | = \sqrt{2i+1}$, $i \in \bbN_0$, it follows that \R{mu_cond} holds with
\be{
\label{LU_ui}
u_i = \sup_{t \in (-1,1)^d} | \phi_i(t) | = \prod^{d}_{j=1} \sqrt{2 i_j+1}.
}

\vspace{1pc}\noindent  \textit{Tensor Legendre polynomials, random sampling from the Chebyshev measure.} In this case, $\nu(t) = 2^{-d}$ and $\mu(t) = \prod^{d}_{j=1} \frac{1}{\pi (1-t^2_j)^{1/2}}$.  Hence
\bes{
u_i = \max \left \{ 1 , \sup_{t \in (-1,1)^d} (\pi/2)^{d/2} \prod^{d}_{j=1} (1-t^2_j)^{1/4} | \phi_{i_j}(t_j) | \right \}.
}
It is known that univariate Legendre polynomials satisfy $(1-t^2)^{1/4} | \phi_0(t) | \leq 1$ and
\be{
\label{1DLeg_envelope}
| \phi_i(t) | (1-t^2)^{1/4} < 2 / \sqrt{\pi},\quad t \in [-1,1],\qquad i \in \bbN,
} 
(see Remark \ref{r:LegWindowBound}).  Therefore we have
\be{
\label{u_i_LC_bound}
u_i \leq (\pi/2)^{d/2} (2/\sqrt{\pi})^{| i|_0},\qquad i \in \bbN^d_0.
}
Unlike the previous cases this bound is not exact for each $i$.  However, it is known that the constant in \R{1DLeg_envelope} cannot be improved.  Hence the bound \R{u_i_LC_bound} is sharp when taken over all $i \in \bbN^d_0$.

\rem{
\label{r:LegWindowBound}
Let $P_{0},P_1,\ldots$ be the classical univariate Legendre polynomials, i.e.\ with normalization $P_i(1) = 1$, so that $\phi_i(t) = \sqrt{2i+1} P_i(t)$.  These polynomials satisfy the following inequality \cite{AntonovHolshevnikovBernstein} (see also \cite{LorchBernstein} and \cite{GautschiBernsteinSharp})
\be{
\label{ClassicalLegEnvelope}
(\sin \theta)^{1/2} |P_i(\cos(\theta)) | < \left (\frac{2}{\pi} \right )^{1/2}(i+1/2)^{-1/2},\quad 0 \leq \theta\leq \pi.
}
Substituting $\phi_i$ and writing $t = \cos(\theta)$ immediately gives \R{1DLeg_envelope}.  Note that \R{ClassicalLegEnvelope} was first proved by Bernstein with the factor $i^{-1/2}$ instead of $(i+1/2)^{-1/2}$ on the right-hand side (see \cite[Thm.\ 7.3.3]{SzegoOrthPolys}).  This weaker inequality has been used several times in the analysis of CS for function approximation \cite{HamptonDoostanCSPCE,Rauhut}.  The sharper bound \R{ClassicalLegEnvelope} allows one to obtain stronger results in the high-dimensional setting for the case of Legendre polynomials with sampling from the Chebyshev measure.  See Corollary \ref{c:LC}.
}

\vspace{1pc} \noindent
In order to formulate \R{fin_min}, we also need to choose the finite index sets $I_K$.  We shall consider the following three standard constructions.  First, the tensor product index set
\bes{
I^{TP}_K = \left \{ i \in \bbN^d_0 : |i|_{\infty} \leq K \right \},
}
where $|i|_{\infty} = \max \{ i_1,\ldots,i_d\}$.  Note that $|I^{TP}_K| = (K+1)^d$.   Although this indexing is arguably the simplest, for moderate $d$ the cardinality of $I^{TP}_K$ is often too large for computations.  A common alternative is the so-called total degree space
\bes{
I^{TD}_{K} = \left \{ i \in \bbN^d_0 : |i|_1 \leq K \right \},
}
where $|i|_1 = i_1+\ldots+i_d$ \cite{NarayanZhouCCP}.  Note that $| I^{TD}_{K} | = \left ( \begin{array}{c} K+d \\ d \end{array} \right )$.  We shall also consider the (isotropic) hyperbolic cross space
\bes{
I^{HC}_K = \left \{ i \in \bbN^d_0 : \prod^{d}_{j=1} (i_j+1) \leq K \right \}.
}
The exact cardinality of this space is harder to quantify, but it is known to satisfy the upper bound
\be{
\label{HCcardinality}
| I^{HC}_K | \leq \min \left \{ 2 K^3 4^d , \E^2 K^{2+\log_2(d)} \right \}.
}
The first inequality is due to \cite[Thm.\ 3.7]{ChernovDungHCCardinality} (using parameters $T=K$, $s=d$, $a=1$ and $\delta = 1/2$), and the second follows from the proof of Theorem 4.9 in \cite{KuhnEtAlApproxMixed}.  See also \cite{ChkifaDownwardsCS}.

At this point, we stress that our main result (Theorem \ref{t:main}) is general, and applies to arbitrary function systems.  We consider Legendre and Chebyshev polynomials with the above samplings since they are popular examples in the literature.  But other cases could also be considered within our framework; for example, Jacobi polynomials.  Our framework and theoretical guarantees also allow for nonpolynomial systems, e.g.\ spherical harmonics or piecewise polynomials.

\section{Main results}\label{ss:main_res}
We now present our main results.  For this, we now introduce the notation $A \gtrsim B$ or $A\lesssim B$ to mean that there exists a constant $C > 0$ independent of all relevant parameters such that $A \geq C B$ or $A \leq C B$ respectively.  In particular, the constant $C$ is independent of the weights used.

\subsection{General recovery guarantee}
To state our results we require two quantities.  First, for weights $w = \{ w_i \}_{i \in I}$ and $u = \{ u_i \}_{i \in I}$ and a finite set $\Delta \subseteq I$, we define
\be{
\label{measurement_condition}
\cM(\Delta ; u , w ) = | \Delta |_u + \max_{i \in I_K \backslash \Delta } \{ u^2_i / w^2_i \} \max \{ | \Delta |_w,1\}.
}
This quantity will play a crucial role in our estimates for the number of measurements required.  Second, for weights $w = \{ w_i \}_{i \in I}$ and $x \in \ell^1_w(I)$, we define
\be{
\label{TK_def}
T_{K,w}(x) = \min \left \{ \| x - \bar{x} \|_{1,w} : \bar{x} \in P_K(\ell^1_w(I)), \| U P_K \bar{x} - y \| \leq \eta \right \}.
}
Loosely speaking, this term determines the additional error incurred due to truncation; that is, in solving the computable minimization problem \R{fin_min} rather than \R{inf_min}.

\thm{
\label{t:main}
Let $K \in \bbN$, $0 < \epsilon < \E^{-1}$, $w = \{ w_i \}_{i \in I}$ be weights, $x \in \ell^1_w(I)$ and $\Delta \subseteq I_K$, $\Delta \neq \emptyset$, be any set with $\min_{i \in I_K \backslash \Delta} \{ w_i \} \geq 1$.  Let $t_1,\ldots,t_m$ be drawn independently from the measure $\mu$.  Then, for all minimizers $\hat{x}$ of \R{fin_min} we have
\be{
\label{main_err}
\| x - \hat{x} \| \lesssim \lambda \sqrt{|\Delta|_w} \left(  \eta / \sqrt{m} + \| x - P_K x \|_{1,u} \right ) + \| x - P_{\Delta} x \|_{1,w} + T_{K,w}(x),
}
with probability at least $1-\epsilon$, provided
\be{
\label{main_est}
m \gtrsim \cM(\Delta ; u , w )  \cdot \log(\epsilon^{-1}) \cdot \log \left (2N \max \left \{ \sqrt{| \Delta |_w }, 1 \right \} \right ),
}
where $N = | I_K|$, $u = \{ u_i \}_{i \in I}$ and $\cM(\Delta ; u,w)$ are as in \R{u_i_def} and \R{measurement_condition} respectively, and $\lambda = 1 + \frac{\sqrt{\log(\epsilon^{-1})}}{\log\left (2N\sqrt{\max \{ | \Delta |_w , 1 \}}\right )}$.
}
We discuss the consequences of this theorem in detail in the following subsections.  However, let us first note that \R{main_err} differs from more standard error bounds in CS in two respects, both of which arise from the truncation.  The first is the term $\| x - P_K x \|_{1,u}$ and the second is the truncation term $T_{K,w}(x)$.  Whilst the first term is small for all sufficiently large $K$, the second term needs to be estimated.  In fact, and much related to this issue, it cannot be taken for granted that \R{fin_min} even has a solution for all $K$, since $U x$ may not lie in the range of $U P_K$.  This is only true in general if $x$ is supported on $I_K$, which is typically not the case in practice.  

Fortunately, both issues can be readily addressed:

\thm{[\cite{Adcockl1pointwise}]
\label{t:trunc_err}
For all sufficiently large $K$, we have $\mathrm{Ran}(U) = \mathrm{Ran}(UP_K)$.  In particular, \R{fin_min} has a solution for all large $K$.
Moreover, suppose that $\mathrm{rank}(U) = r \leq m$ and $K$ is sufficiently large so that $\mathrm{rank}(UP_K) = r$.  If $x \in \ell^1_{\tilde{w}}(I)$ then
\bes{
T_{K,w}(x) \leq \left ( 1 +  \| P_K w \| / \sigma_r \right ) \| x - P_K x \|_{1,w},
}
where $\sigma_{r}$ is the $r^{\rth}$ singular value of $U P_K$.
}

These results imply the following.  Provided $K$ is chosen sufficiently large such that $1/\sigma_{r}$ is small and finite, then the truncated problem \R{fin_min} not only has a solution but the additional error $T_{K,w}(x)$ due to truncation is bounded by the tail error $\| x - P_K x \|_{1,w}$ multiplied by $\| P_K w \|$.  Crucially, the condition $1/\sigma_{r} < \infty$ is independent of $x$ (and therefore $f$) and depends only on the data sample points $\{ t_i \}^{m}_{i=1}$ and the system $\{ \phi_i \}_{i \in I}$.  It can also be easily checked numerically.  For theoretical guarantees relating $K$ to the number of samples $m$ to ensure this condition, we refer to \cite{Adcockl1pointwise}.

\rem{
Note that the term $ \| P_K w \| \| x - P_K x \|_{1,w} \rightarrow 0$ as $K \rightarrow \infty$ given some additional summability of the coefficients $x$.  For example, suppose that $I = \bbN$, $I_K = \{1,\ldots,K\}$ and the weights $w_i$ are nondecreasing.  Then it is straightforward to see that $\| P_K w \| \| x - P_K x \|_{1,w} \leq \| x - P_K x \|_{1,\tilde{w}}$, where $\tilde{w}_i = \sqrt{i} w^2_i$.  Hence this term tends to zero as $K \rightarrow \infty$ provided $x \in \ell^1_{\tilde{w}}(\bbN)$.
}

\rem{
\label{r:Myestimateisbetter}
In the language of CS, Theorem \ref{t:main} in an example of a nonuniform recovery guarantee \cite{FoucartRauhutCSbook}.  For uniform guarantees, we refer to \cite{ChkifaDownwardsCS,RauhutWardWeighted}.  Note that the guarantees in \cite{ChkifaDownwardsCS,RauhutWardWeighted} allow the error to be estimated in the stronger norm $\ell^1_w$-norm, whereas in \R{main_err} uses the weaker $\ell^2$-norm.  This is a standard discrepancy between uniform and nonuniform-based CS analyses \cite{FoucartRauhutCSbook}.  On the other hand, nonuniform recovery arguments are more flexible, in the sense that they can be used to derive recovery guarantees for arbitrary support sets $\Delta$ without specifying a sparsity model (e.g.\ sparsity or weighted sparsity).  See \cite{BigotBlockCS,AdcockChunParallel} for related work in this direction.  It is this flexibility that allows us to derive the bound \R{main_est}.  As is also typical, our nonuniform guarantee \R{main_est} involves fewer log factors than corresponding uniform guarantees.
}

\rem{
\label{r:Myconditionsaremoregeneral}
The considerations of the previous remark aside, the conditions of Theorem \ref{t:main} are also more general than those of \cite{ChkifaDownwardsCS,RauhutWardWeighted}.  Specifically, we do not require the weights to satisfy $w_i \geq u_i$ \cite{RauhutWardWeighted} or $w_i = u_i$ \cite{ChkifaDownwardsCS}, and we do not impose \textit{a priori} estimates on the expansion tail (recall the discussion in \S \ref{ss:current_approaches} and \S \ref{ss:contributions}).
}

\subsection{A weighted sparsity recovery guarantee}

In the absence of noise, Theorem \ref{t:main} states that $x$ is recovered up to an error proportional $\| x - P_{\Delta} x \|_{1,w}$, i.e.\ the norm of the coefficients of $x$ lying outside $\Delta$, provided the number of measurements is, up to log factors, proportional to
\bes{
\cM(\Delta;u,w) =  | \Delta |_u + \max_{i \in I_K \backslash \Delta} \{ u^2_i / w^2_i \} \max \{ | \Delta |_w,1 \}.
}
This bound may at first sight appear rather obscure, since it does not depend solely on the (weighted) sparsity of $x$.  However, the generality of this bound will be useful in subsequent sections to analyze the performance of weighted $\ell^1$ minimization in different scenarios (recall \S \ref{ss:contributions}).  First, though, we note that Theorem \ref{t:main} immediately implies a weighted sparsity recovery guarantee, similar to that introduced in \cite{RauhutWardWeighted}:

\cor{
\label{c:main_inf_nonlinear}
Let $w = \{ w_i \}_{i \in I}$ be weights satisfying $w_i \geq u_i$, $\forall i \in I$, where $u = \{ u_i \}_{i \in I}$ is as in \R{u_i_def} with $\mu = \nu$, i.e.\ $u_i = \| \phi_i \|_{L^\infty}$.  Let $0 < \epsilon < \E^{-1}$, $K \in \bbN$, $x \in \ell^1_w(I)$ and suppose that $t_1,\ldots,t_m$ are drawn independently  from the measure $\mu$.  Suppose that
\be{
\label{main_est_inf_nonlinear}
m \gtrsim s \cdot \log(\epsilon^{-1}) \cdot \log(2N \sqrt{s}),
}
where $N = | I_K|$.  Then, for all minimizers $\hat{x}$ of \R{fin_min}, we have
\bes{
\| x - \hat{x} \| \lesssim \lambda \sqrt{s} \left ( \eta / \sqrt{m} + \| x - P_K x \|_{1,u} \right ) + \sigma_{s,K}(x)_{1,w} + T_{K,w}(x),
}
with probability at least $1-\epsilon$, where $\lambda = 1 + \frac{\sqrt{\log(\epsilon^{-1})}}{\log \left ( 2 N \sqrt{s} \right ) }$, $T_{K,w}(x)$ is as in \R{TK_def} and
\be{
\label{best_s_err}
\sigma_{s,K}(x)_{1,w} = \min \left \{ \| x - P_{\Delta} x \|_{1,w} : \Delta \subseteq I_K, | \Delta |_{w} \leq s \right \},
}
is the best weighted $s$-sparse approximation error.
}

\prf{
Let $\Delta \subseteq I_K$, $|\Delta|_{w} \leq s$ be such that $\| x - P_{\Delta} x \|_{1,w} = \sigma_{s,K}(x)_{1,w}$.  Since $w_i \geq u_i$ we have $| \Delta |_u \leq | \Delta |_w$ and $\cM(\Delta ; u,w) \leq 2 | \Delta |_w$.  The result now follows immediately from Theorem \ref{t:main}.
}

This result shows that \R{fin_min} attains the best weighted nonlinear approximation error $\sigma_{s,K}(x)_{1,w}$, up to a constant, using a number of measurements $m$ scaling linearly with $s$.  Note that Corollary \ref{c:main_inf_nonlinear} is similar to results found in \cite{RauhutWardWeighted}, except for the differences mentioned in Remarks \ref{r:Myestimateisbetter} and \ref{r:Myconditionsaremoregeneral}.

\rem{
The reader will have noticed that \R{best_s_err} is not the true best weighted $s$-sparse approximation error, but rather the best weighted $s$-sparse approximation error up to some finite range $I_K$.  Since coefficients outside $I_K$ do not form part of the optimization problem \R{fin_min}, this definition is natural (in fact, it has been implicitly assumed in all prior theoretical analysis of CS for function approximation).  More fundamentally, one cannot expect to stably recover arbitrary $s$-sparse vectors whose coefficients can range over the whole of $I$ (a countable index set) when taking only a finite number of samples.  This lack of \textit{instance optimality} in infinite dimensions is discussed in \cite{BAACHGSCS,BourrierEtAlInstance}.  

Having said this, suppose that the weights $w_{i}$ are increasing in the sense that $\min_{i \in I \backslash I_K} \{ w_i \} \rightarrow \infty$ as $K \rightarrow \infty$.  Then this condition can be removed since sets of weighted cardinality $| \Delta |_w \leq s$ cannot have arbitrary range.  In particular, if $K$ is chosen so that $\min_{i \in I \backslash I_K} \{ w_i \} \geq \sqrt{s}$ then one can replace $\sigma_{s,K}(x)_{1,w}$ with the true best weighted $s$-sparse approximation error
\bes{
\sigma_{s}(x)_{1,w} = \min \left \{ \| x - P_{\Delta } \|_{1,w} : \Delta \subseteq I,\ | \Delta |_w \leq s \right \}.
}
}

\section{Consequences of Theorem \ref{t:main}}\label{s:consequences}
In this section, we discuss the main consequences of Theorem \ref{t:main} as listed in \S \ref{ss:contributions}.

\subsection{Recovery guarantees for linear models}\label{ss:linear_models}
As discussed in \S \ref{ss:recov_intro}, it may be the case that the most significant coefficients of a given function are located at the lowest indices with respect to some ordering.  Suppose for simplicity that $\Delta = \{ j : x_j \neq 0 \} = I_M$ for some $M$.  Then, according to the weighted sparsity recovery guarantee (Corollary \ref{c:main_inf_nonlinear}), the number of measurements needed to recover $x$ is, up to log factors, proportional to
\bes{
 | I_M |_{w} = \sum_{i \in I_M} w^2_i.
}
However, this condition depends on the optimization weights $w$ and deteriorates as they increase.  Conversely, the following is a straightforward  corollary of Theorem \ref{t:main} gives a sharper estimate:

\cor{
\label{c:linear_sharp}
Let $w = \{ w_i \}_{i \in I}$ be any weights satisfying
\be{
\label{weights_cond_2}
\max_{i \in I_R} \{ w_i / u_i \} \asymp \inf_{i \in I \backslash I_R } \{ w_i / u_i \} \asymp R^{\upsilon},\qquad R \rightarrow \infty,
}
for some $\upsilon \geq 0$, where $u = \{ u_i \}_{i \in I}$ is as in \R{u_i_def}.  Let $0 < \epsilon < \E^{-1}$, $K \in \bbN$, $x \in \ell^1_w(I)$ and suppose that $t_1,\ldots,t_m$ are drawn independently from the measure $\mu$.  If $I_M \subseteq I_K$, $\min_{i \in I_K \backslash I_M } \{w_i \} \geq 1$ and
\be{
\label{main_est_linear}
m \gtrsim | I_M |_{u} \cdot \log(\epsilon^{-1}) \cdot \left ( \log\left (2N \sqrt{| I_M |_u} \right ) + \nu \log(M) \right ),
}
then, for any minimizer $\hat{x}$ of \R{fin_min}, we have
\be{
\label{main_err_inf}
\| x - \hat{x} \| \lesssim \lambda \sqrt{|I_M |_w} \left ( \eta /\sqrt{m} + \| x- P_K x \|_{1,u} \right )  + \| x - P_{M} x \|_{1,w} + T_{K,w}(x),
}
with probability at least $1-\epsilon$, where $\lambda = 1 + \frac{\sqrt{\log(\epsilon^{-1})}}{\log(2N \sqrt{| I_M |_w})}$.
}
\prf{
We use Theorem \ref{t:main} with $\Delta = I_M$.  Observe that $|I_M|_{w} \lesssim M^{2 \nu} | I_M |_{u}$ and
\bes{
\max_{i \in I_K \backslash I_M } \{ u^2_i / w^2_i \} \max \left \{ | I_M |_w , 1 \right \} \lesssim M^{-2 \nu} \max \left \{ M^{2 \nu} | I_M |_u , 1 \right \} \lesssim | I_M |_u,
}
where in the final step we use the fact that $u_i \geq 1$, $\forall i$.  The result now follows immediately.
}

Observe that \R{main_est_linear} depends only on the intrinsic weights $u$ and is independent of the optimization weights $w$, provided these weights satisfy \R{weights_cond_2}.  Loosely speaking, this means that the $w_i$'s must grow at least as fast as the $u_i$'s.

\examp{
\label{ex:ex}
Let $d=1$, $I = \bbN$ and $I_M = \{1,\ldots,M\}$ for $M \in \bbN$.  As mentioned in \S \ref{ss:recov_intro}, an oscillatory function (with frequency of oscillation $\ord{M}$) typically has $x_i = \ord{1}$ for $i = 1,\ldots,M$ and $x_i \approx 0$ for $i > M$.  Hence a good approximation of such a function occurs only if the first $M$ coefficients are accurately recovered.
Suppose now that the weights $w_i = i^{\alpha}$ for some $\alpha > 0$.  Then according to Corollary \ref{c:main_inf_nonlinear} the number of measurements needed is roughly $M^{2 \alpha +1}$.  Thus, more measurements are apparently required for more rapidly growing weights, at odds with the results shown in Fig.\  \ref{f:ErrIntroduction}.  Suppose now that the intrinsic weights $u_{i} \asymp i^{\beta}$ as $i \rightarrow \infty$ for some $\beta \geq 0$.  If $\alpha \geq \beta$ then Corollary \R{c:linear_sharp} (with $\nu = \alpha - \beta$) gives that the number of measurements is proportional (up to log factors) to $|I_M|_{u} \asymp M^{2 \beta + 1}$, regardless of $\alpha$.  

In the case of univariate Chebyshev or Legendre polynomials with sampling from the Chebyshev measure -- in which case the $u_i$'s are uniformly bounded (see \S \ref{ss:main_examp}) and therefore $\beta = 0$  -- this result gives a linear scaling of $m$ with $M$, up to log factors, regardless of the choice of $\alpha$.  Similarly, for Legendre polynomials with sampling from the uniform measure (in which case $\beta = 1/2$), one deduces a quadratic scaling of $m$ with $M$.  Both results are essentially optimal; see Remark \ref{r:quadscalingsharp}.
}

Corollary \ref{c:linear_sharp} can also be used to assert similar results in higher dimensions.  For example, suppose that the coefficients  $x_i$ satisfy
\bes{
|x_i| \asymp  \prod^{d}_{j=1} (i_j+1)^{-\beta},
}
for some $\beta > 0$, as is reasonable in some cases \cite{WebsterPolyApproxi}.  Then the significant coefficients lie in a hyperbolic cross $I_M = I^{HC}_M$.  Suppose now that the weights $w_i$ are chosen as $w_i = u_i \prod^{d}_{j=1} (i_j+1)^\alpha$.  Then \R{weights_cond_2} holds with $\nu = \alpha$, leading to a measurement condition proportional to $|I^{HC}_M|_u$.  As in Example \ref{ex:ex}, this is independent of the parameter $\alpha$.  Similarly, if the coefficients $|x_i| \asymp \rho^{-|i|_1}$ for some $\rho > 1$ then one may take $I_M = I^{TD}_{M}$ to be a total degree index set.  If the weights are chosen as $w_i = u_i ( |i|_1 +1 )^\alpha$ then Corollary \ref{c:linear_sharp} gives a measurement condition proportional to $|I^{TD}_M|_u$, which is once more independent of $\alpha$.

\begin{figure}
\begin{center}
\begin{tabular}{ccc}
LU, $d=2$ & CC, $d=3$ & LC, $d=4$ 
\\
\includegraphics[width=5.25cm]{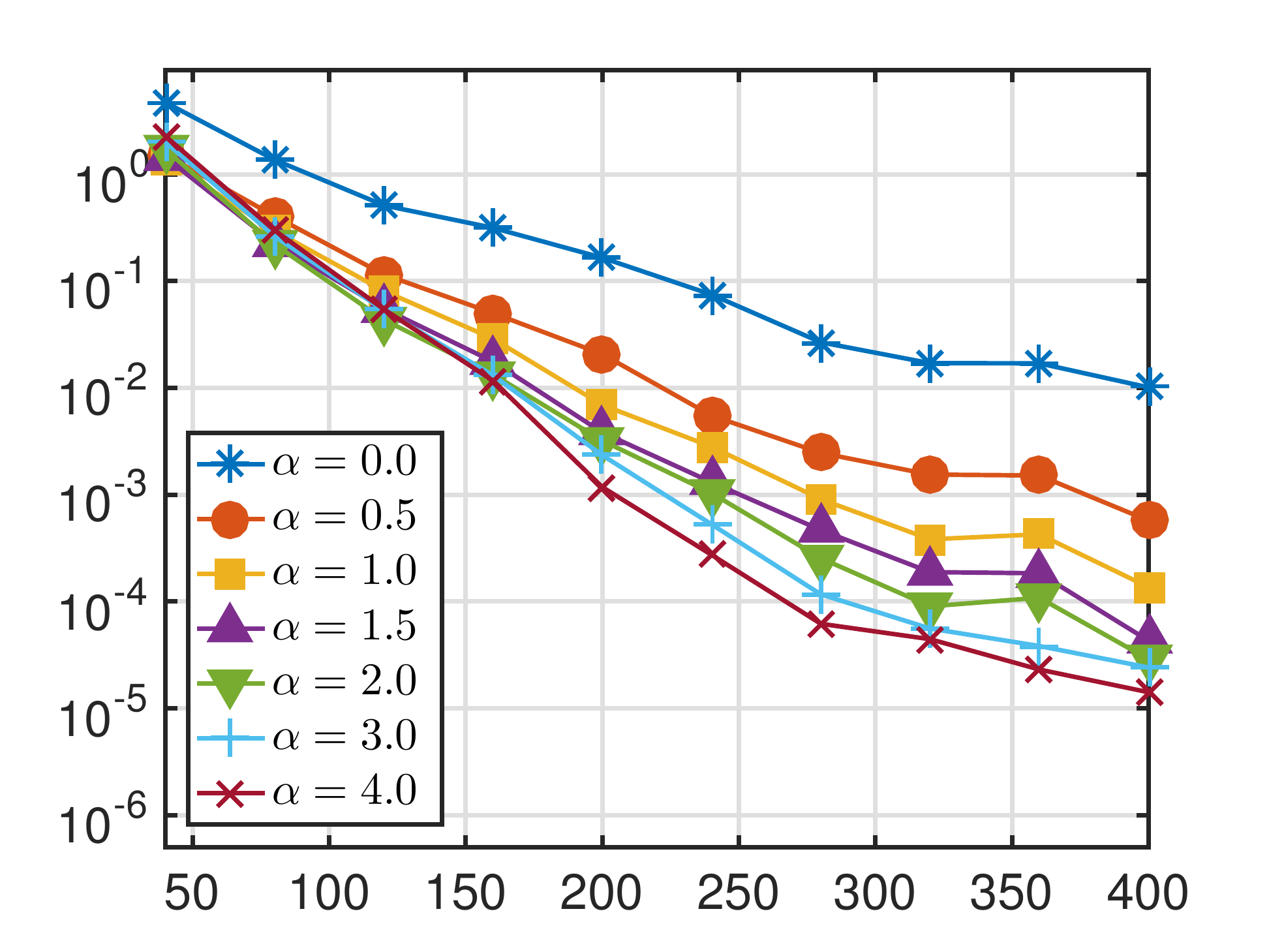} &
\includegraphics[width=5.25cm]{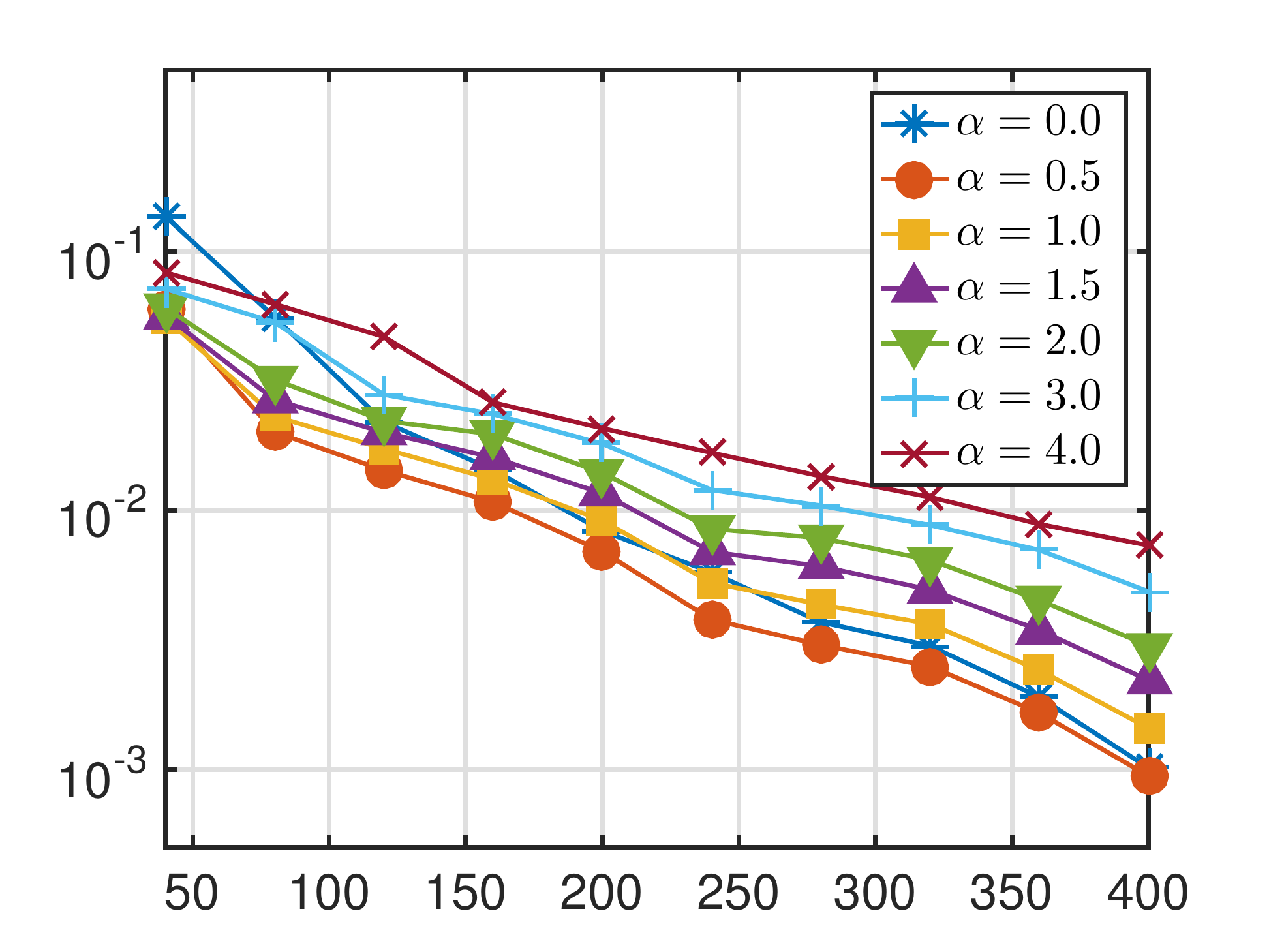} &
\includegraphics[width=5.25cm]{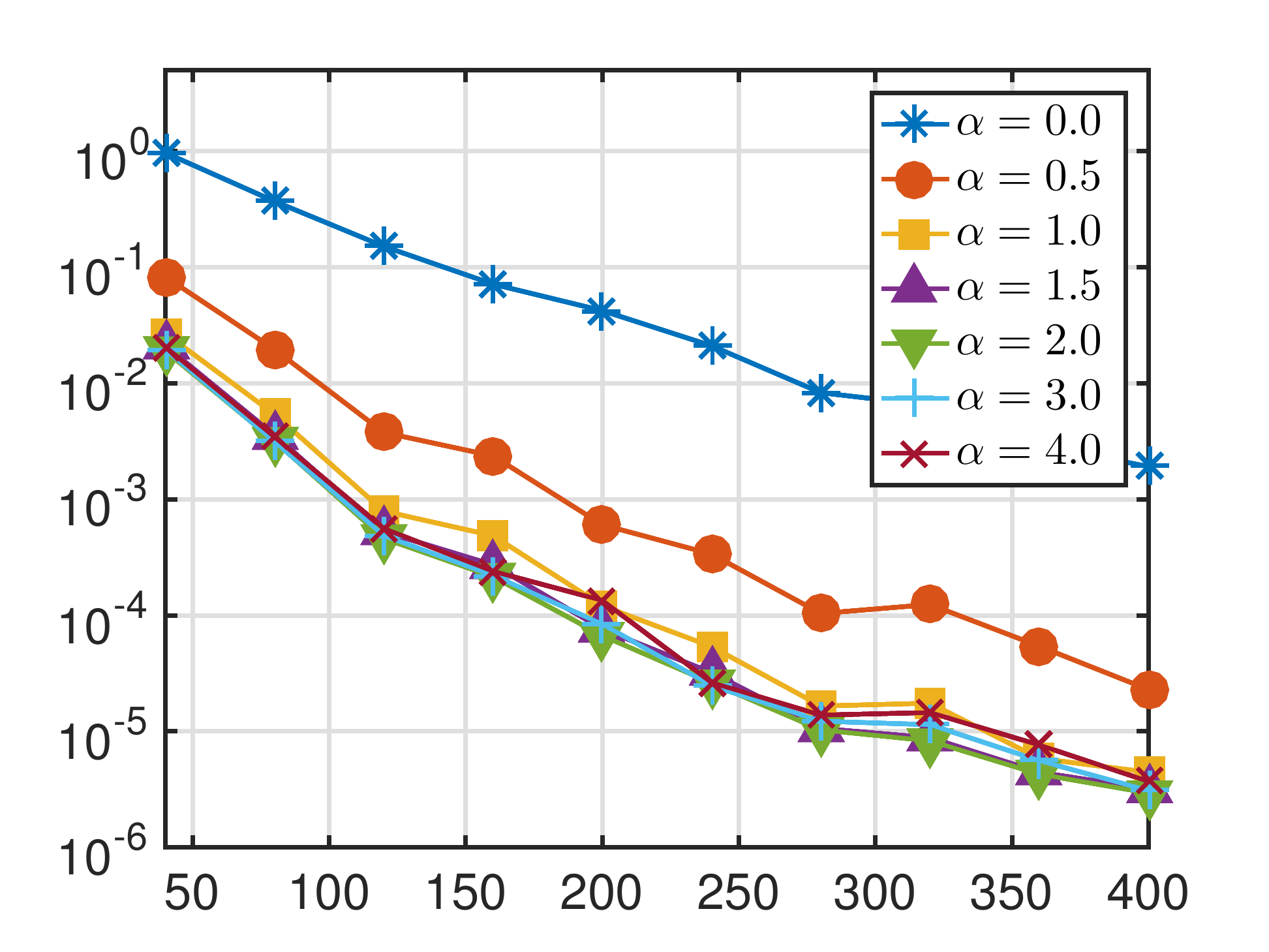} 
\\
\includegraphics[width=5.25cm]{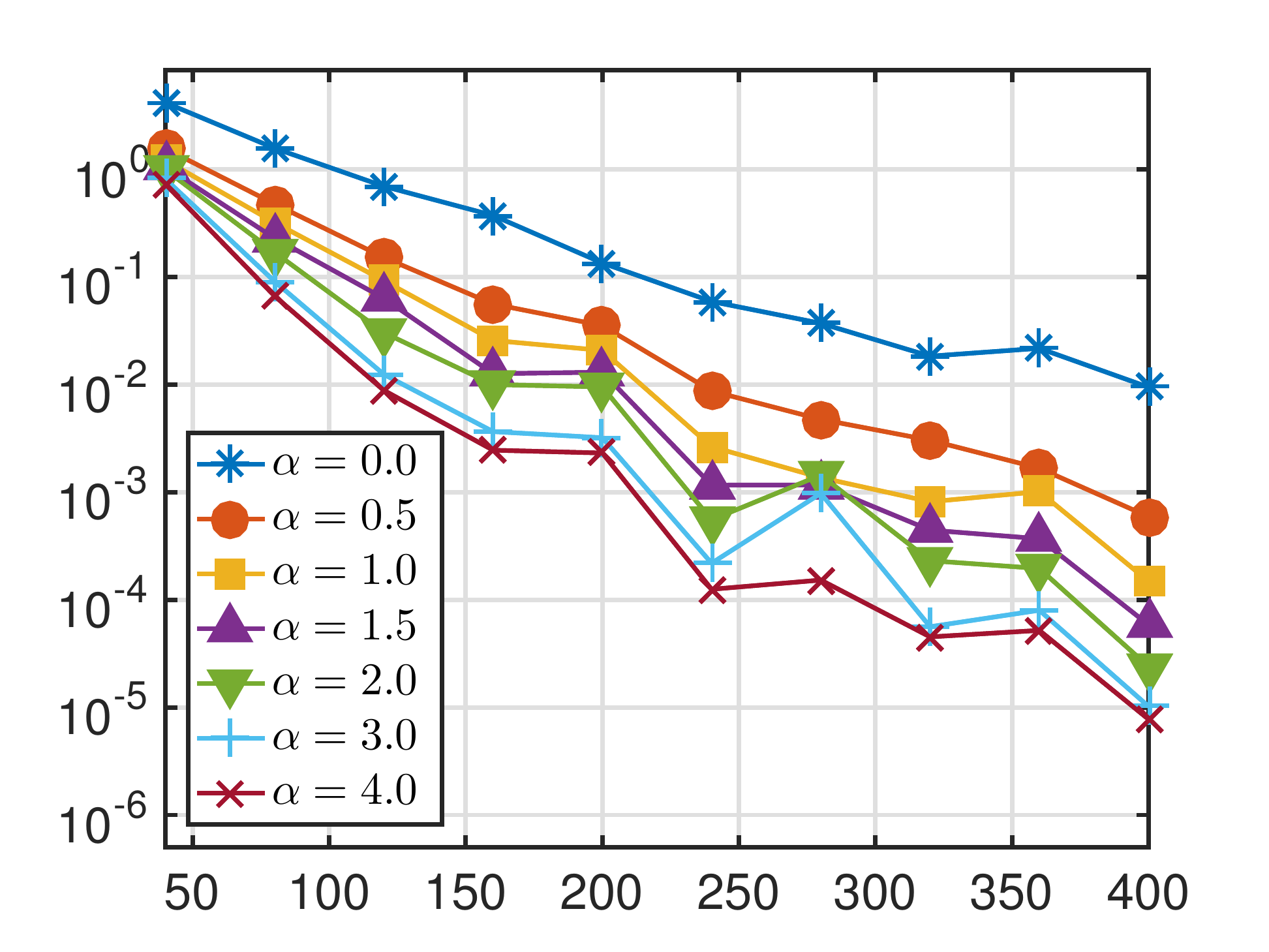} &
\includegraphics[width=5.25cm]{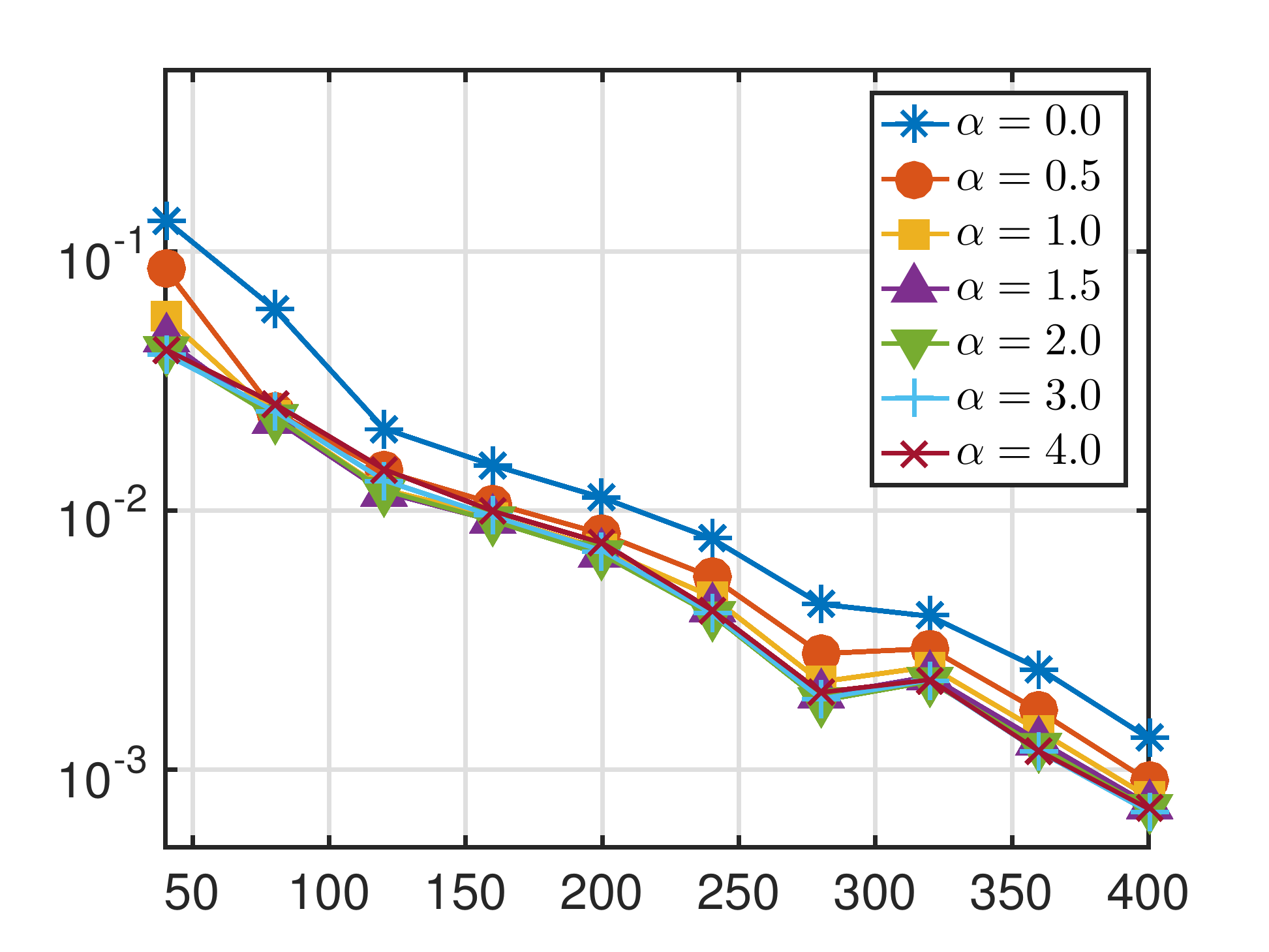} &
\includegraphics[width=5.25cm]{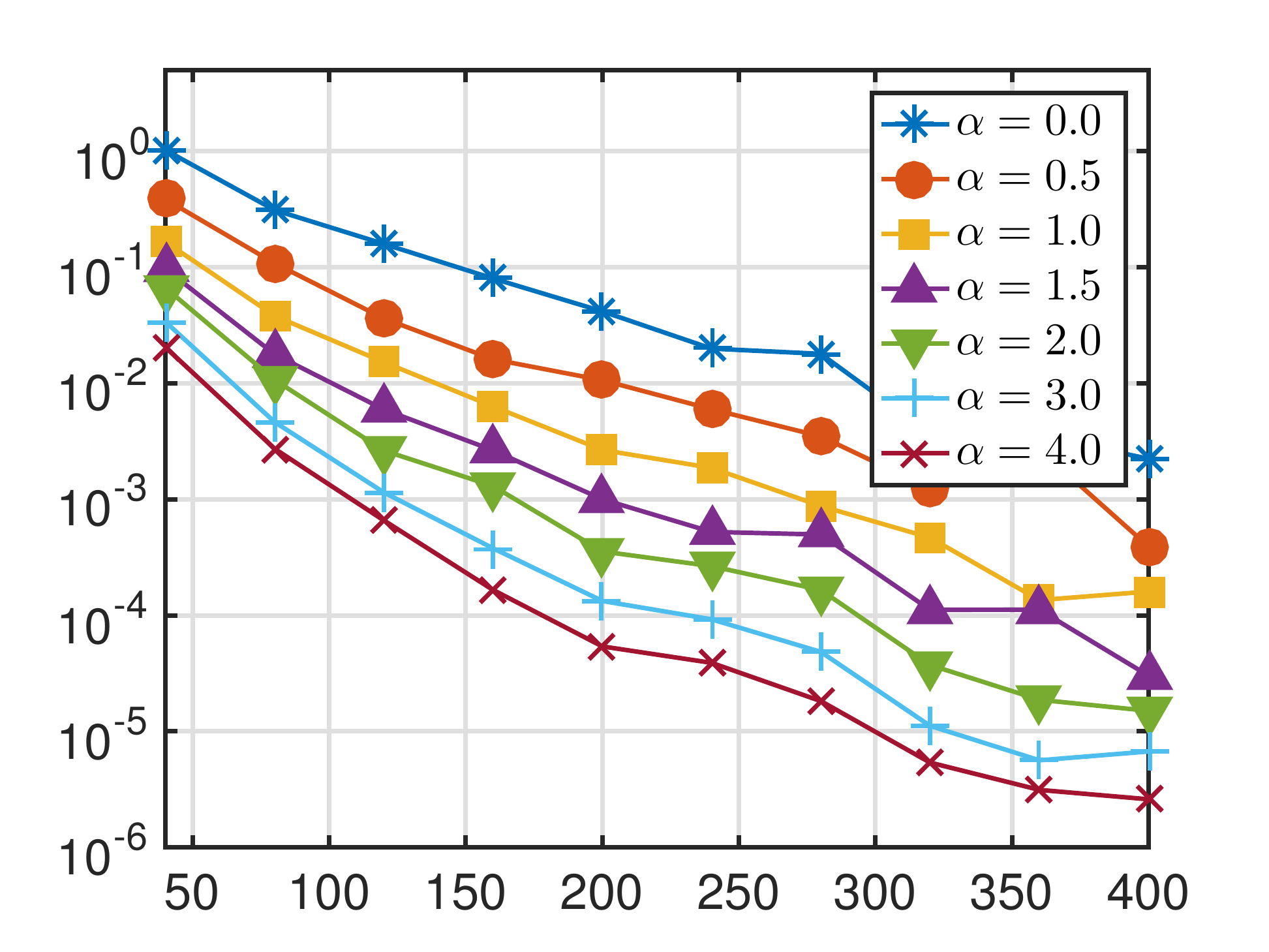} 
\\
$f(t) = \exp(2 t_1) \cos(3 t_2)$  & $f(t) = \sin(0.5\exp(t_1t_2t_3))$ & $f(t) = \exp\left (-\frac{t_1+t_2+t_3+t_4}{6} \right)$
\end{tabular}
\end{center}
\caption{
The error $\| f - \tilde{f} \|_{L^\infty}$ (averaged over $50$ trials) against $m$ for Chebyshev (C) or Legendre (L) polynomials with points drawn from the Chebyshev (C) or uniform (U) measure.  A total degree index set of degree $K$ was used, where $(d,K) = (2,44),(3,17),(4,10)$, and the weights were taken to be either $w_i = \prod^{d}_{j=1} (i_j+1)^{\alpha}$ (top row) or $w_i = (|i|_1 + 1)^{\alpha}$ (bottom row) for various $\alpha \geq 0$.
}
\label{f:dDlinearmodel_weights}
\end{figure}

Numerical illustrations of these results are given in Fig.\ \ref{f:ErrIntroduction} (the univariate case) and Fig.\ \ref{f:dDlinearmodel_weights} (the multivariate case).

\subsection{The choice $1 \leq w_i \leq u_i$}
In this and the following two subsections we turn our attention to using Theorem \ref{t:main} to understand the benefits that weights convey, as opposed to just showing that they lead to no deterioration in the recovery guarantee.   
We commence with the following straightforward result:

\prop{
Let $w_{i} = (u_i)^{\theta}$ for some $0 \leq \theta \leq 1$.  Then
\bes{
| \Delta |_u + \max_{i \in I_K} \{ u^2_i / w^2_i \} | \Delta |_w \leq | \Delta |_u + \max_{i \in I_K } \{ u^2_i \} | \Delta |.
}
}
\prf{
We have $\max_{i \in I_K} \{ u^2_i / w^2_i \} | \Delta |_w \leq \max_{i \in I_K } \{ u^{2(1-\theta)}_i \} \max_{i \in I_K } \{ u^{2 \theta}_i \} | \Delta | = \max_{i \in I_K } \{ u^2_i \} | \Delta |$, as required.
}

This result shows that choosing weights $w_i$ scaling with the $u_i$'s cannot worsen the recovery guarantee (except possibly in the log factor) over that of the unweighted ($\theta = 0$) case.  Of particular interest is the extreme case $\theta = 1$, in which case we have 
\bes{
| \Delta |_u + \max_{i \in I_K} \{ u^2_i / w^2_i \} | \Delta |_w = 2 | \Delta |_u,\qquad w_i = u_i.
}
We expect this to be significantly smaller than the corresponding estimate for the unweighted case
\bes{
| \Delta |_u + \max_{i \in I_K} \{ u^2_i / w^2_i \} | \Delta |_w =| \Delta |_u + \max_{i \in I_K } \{ u^2_i \} | \Delta |,\qquad w_i = 1,
}
whenever support set $\Delta$ does not contain too many high indices, so that the weighted cardinality $| \Delta |_u$ is not too large in comparison to the maximum of the $u_i$'s over the range $I_K$.  In the next subsection we will see several concrete examples of this in the case of polynomial approximation.

\subsection{Recovery guarantees for tensor Chebyshev and Legendre expansions}\label{ss:recovery_examples}
We now specialize our focus to the case of tensor Chebyshev and Legendre polynomial expansions.  Following on from the previous subsection, we consider the cases $w_i = 1$ and $w_i = u_i$ respectively.
Throughout we let $I = \bbN^d_0$ and for the truncated spaces we consider total degree spaces $I_K = I^{TD}_{K}$.  

It will also be useful to first recall the definition of a lower set:

\defn{
\label{d:lower_set}
A set $\Delta \subseteq \bbN^{d}_{0}$ is lower if whenever $i = (i_1,\ldots,i_d) \in \Delta$ and $i' = (i'_1,\ldots,i'_d) \in \bbN^d_0$ satisfies $i'_j \leq i_j$, $j=1,\ldots,d$, then $i' \in \Delta$.
}
Lower (or sometimes referred to as downwards closed) sets are well-known constructions in multivariate polynomial approximation, since in practice that the support sets of polynomial coefficients are often described by such sets.  See \cite{ChkifaEtAl,ChkifaDownwardsCS,CohenDeVoreSchwabFoCM,CohenDeVoreSchwabRegularity,MiglioratiJAT,MiglioratiEtAlFoCM} and references therein.

\vspace{1pc} \noindent
\textit{Tensor Chebyshev polynomials, random sampling from the Chebyshev measure.} 

\cor{
\label{c:CC}
Let $\nu(t) = \mu(t) = \prod^{d}_{j=1} \frac{1}{\pi(1-t^2_j)^{1/2}}$.  Then, for any $\Delta \subseteq I_K$ with $| \Delta | \leq s$ we have 
\be{
\label{CC_w1}
\cM(\Delta ; u , 1 ) \leq 2^{\min\{d,K\}+1} s,
}
provided $I_K = I^{TD}_{K}$ is the total degree index set.  If $\Delta \subseteq I_K$ is also a lower set then
\be{
\label{CC_w2}
 \cM(\Delta ; u , u) \leq 2  s^{\log(3)/\log(2)} , 
}
regardless of the choice of $I_K$.  In other words, lower sets of cardinality $s$ can be recovered via weighted $\ell^1$ minimization with weights $w_i = 2^{|i|_0/2}$ from a number of measurements that is independent of $d$ for large $d$ and proportional to $s^{\log(3)/\log(2)}$. 
}
\prf{
The weights $u_i$ are given by \R{CC_ui} in this case.  Since $\Delta$ is a subset of the total degree space, we have $|i|_0 \leq \min \{ d,K\}$ for $i \in \Delta$.  The first result now follows from \R{CC_ui} and the definition of $\cM$.  For the second we note that $\sum_{i \in \Delta} 2^{|i|_0} \leq |\Delta|^{\log(3)/\log(2)}$ for any lower set \cite{ChkifaEtAl}.
}

\begin{figure}
\begin{center}
\begin{tabular}{ccc}
\includegraphics[width=5.25cm]{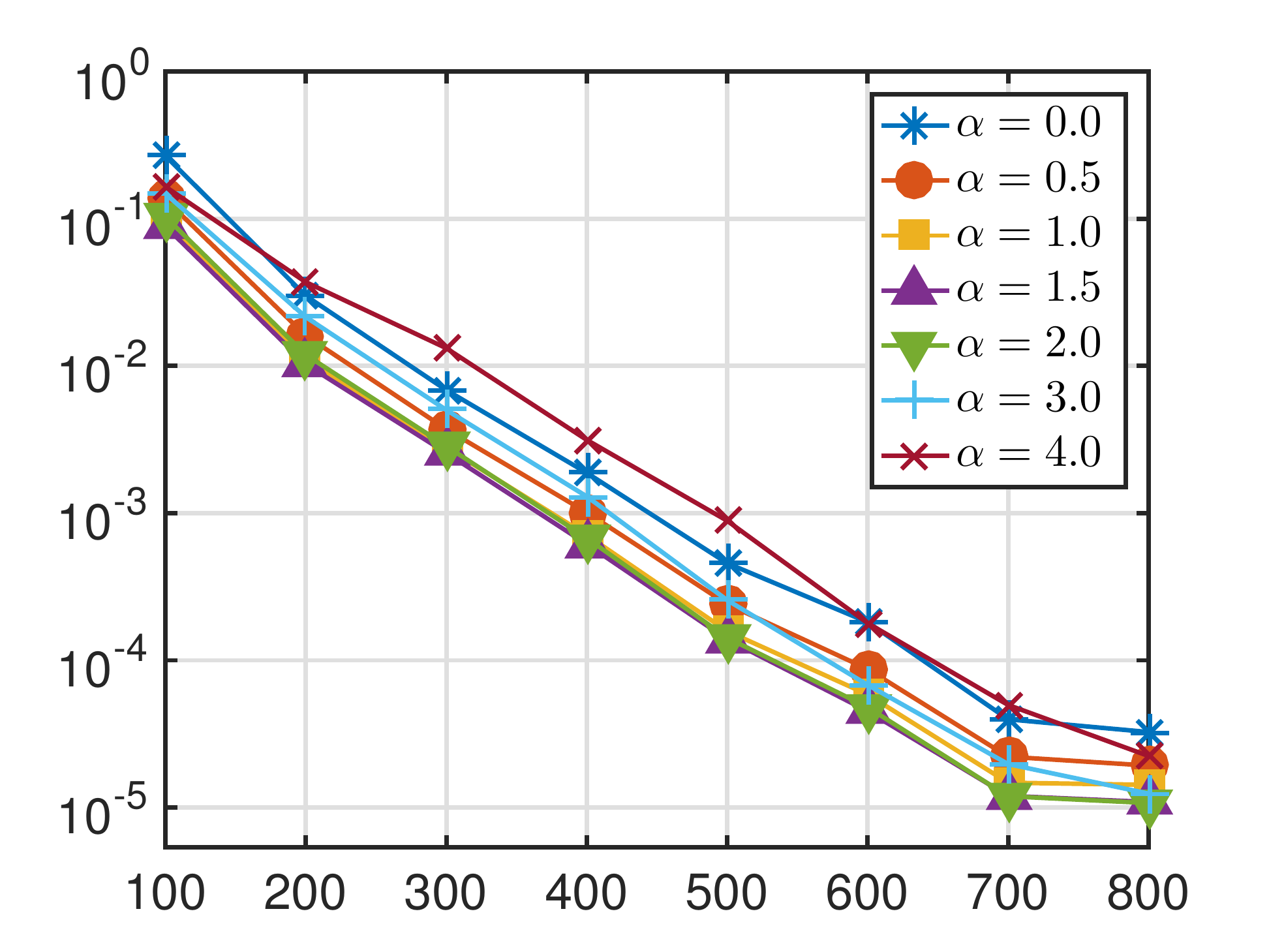} &
\includegraphics[width=5.25cm]{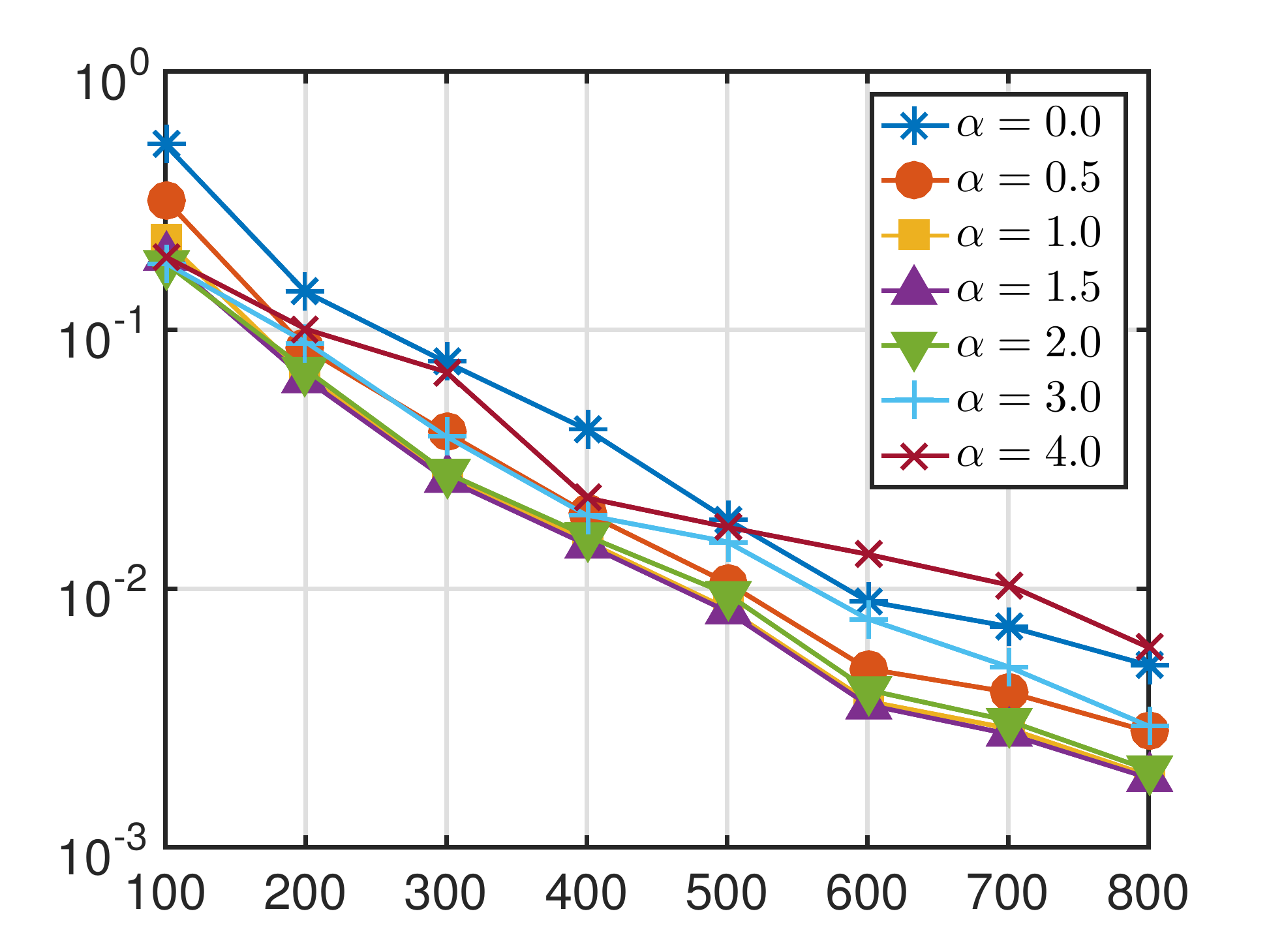} &
\includegraphics[width=5.25cm]{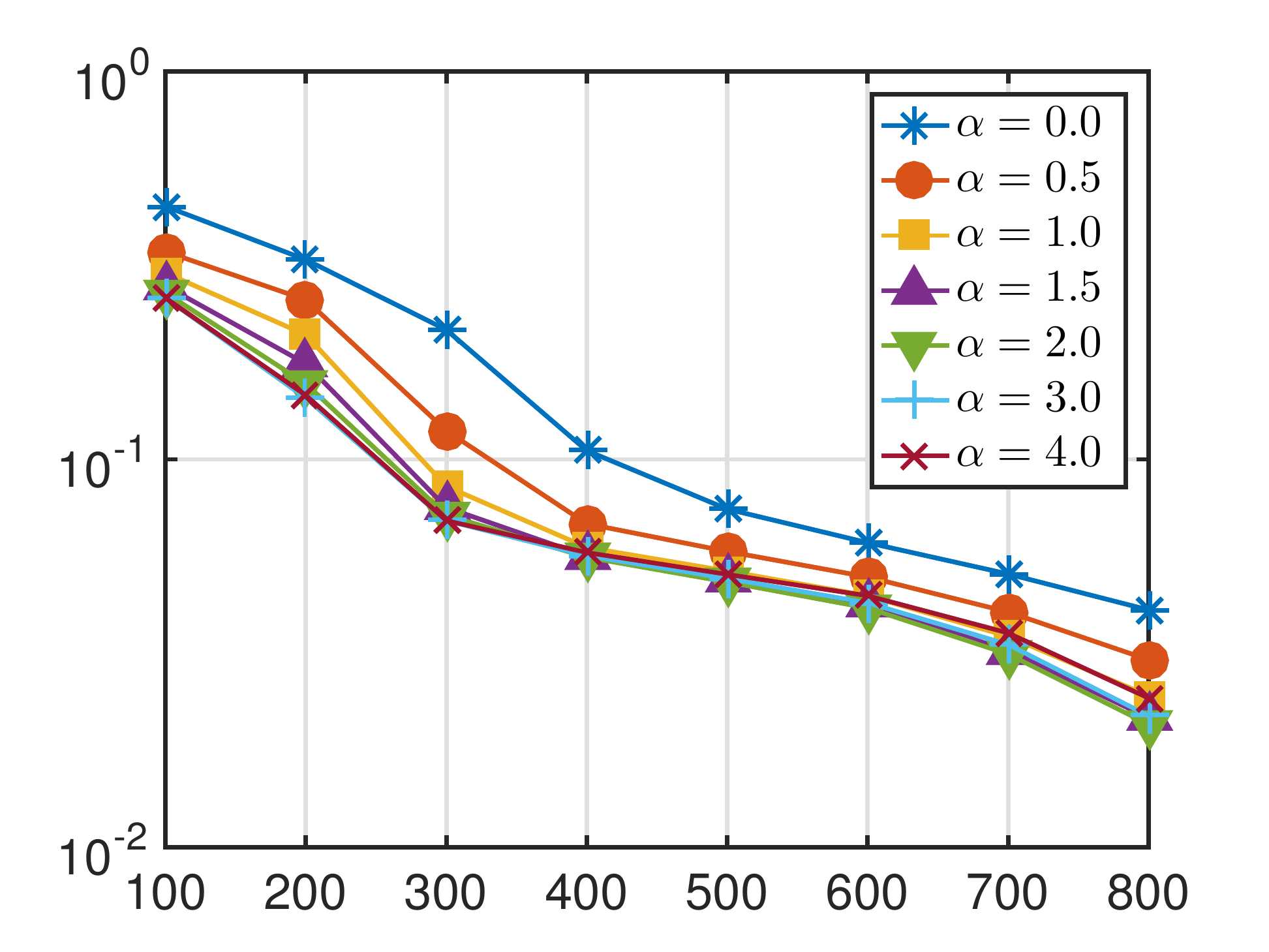}
\\
$d=3$ & $d=5$ & $d=10$
\end{tabular}
\end{center}
\caption{
The error $\| f - \tilde{f} \|_{L^\infty}$ (averaged over $50$ trials) against $m$ for $f(t) = \exp(-(t_1+\ldots+t_d)/(2d))$ and Chebyshev polynomials with points drawn from the Chebyshev measure.  A total degree index set of degree $K$ was used, where $(d,K) = (3,24),(5,10),(10,5)$.  The weights were taken to be $w_i = (u_i)^{\alpha}$ for various $\alpha$, where $u_i$ is as in \R{CC_ui}.
}
\label{f:CC_dimcomp}
\end{figure}

This result demonstrates the advantage of setting the weights $w_i = u_i$.  In high dimensions, one can recover lower sets of coefficients using a number of measurements that is (up to log factors) independent of the dimension.  Note that the scaling $\log(3)/\log(2)$ is sharp in the sense that it agrees with the best known estimates for recovering a fixed, known lower set via discrete least-squares \cite{ChkifaEtAl}.  A numerical illustration of this result is given in Fig.\ \ref{f:CC_dimcomp}.  In all dimensions, setting the weights as $w_i = (u_i)^{\alpha}$ for some $\alpha > 0$ leads to a smaller approximation error, with the choice $w_i = u_i$ giving amongst the smallest.  This is in good agreement with the above corollary.

We remark that \R{CC_w1} was first obtained in \cite{YanGuoXui_l1UQ} (see also \cite{HamptonDoostanCSPCE}) and \R{CC_w2} has also been presented in \cite{ChkifaDownwardsCS}.  Although our main condition is the same, our analysis improves on both results by removing the requirement for \textit{a priori} tail estimates (recall \S \ref{ss:current_approaches}).  As discussed, our recovery guarantees also exhibit fewer log factors (see Remark \ref{r:Myestimateisbetter}).

\vspace{1pc} \noindent
\textit{Tensor Legendre polynomials, random sampling from the uniform measure.}

\cor{
\label{c:LU}
Let $\nu(t) = \mu(t) = 2^{-d}$.  Then, for any $\Delta \subseteq I_K$ with $| \Delta | \leq s$ we have 
\be{
\label{LU_w1}
\cM(\Delta ; u,1) \leq 2 \times 3^K s,
}
provided $I_K = I^{TD}_{K}$ is the total degree space of degree $K$.  Conversely, if $\Delta \subseteq I_K$ is a lower set then
\be{
\label{LU_w2}
\cM(\Delta ; u,u) \leq 2 s^2,
}
regardless of the choice of $I_K$.
}
\prf{
Note that the weights $u_i$ satisfy \R{LU_ui}.  For \R{LU_w1} we recall from \cite{YanGuoXui_l1UQ} that $u^2_i \leq 3^K$, $i \in I_K$, whenever $I_K$ is the total degree space.  For \R{LU_w2}, we recall that $| \Delta |_u \leq s^2$ for lower sets \cite{ChkifaEtAl}.
}

As in the previous case, this result clearly illustrates the benefits of choosing weights $w_i = u_i$.  Note that \R{LU_w1} was first obtained in \cite{YanGuoXui_l1UQ} (see also \cite{HamptonDoostanCSPCE}) and \R{LU_w2} has also been given in \cite{ChkifaDownwardsCS}.  We remark in passing that the estimate \R{LU_w1} is sharp when $d > K$, but ceases to be sharp when $d \leq K$.  For a better bound in this regime, see \cite{HamptonDoostanCSPCE}.  Also as in the previous setting, we note that the bound \R{LU_w2} is sharp in the sense it agrees with the best known estimates for recovery of a fixed lower set via discrete least squares \cite{ChkifaEtAl}.

\begin{figure}
\begin{center}
\begin{tabular}{ccc}
\includegraphics[width=5.25cm]{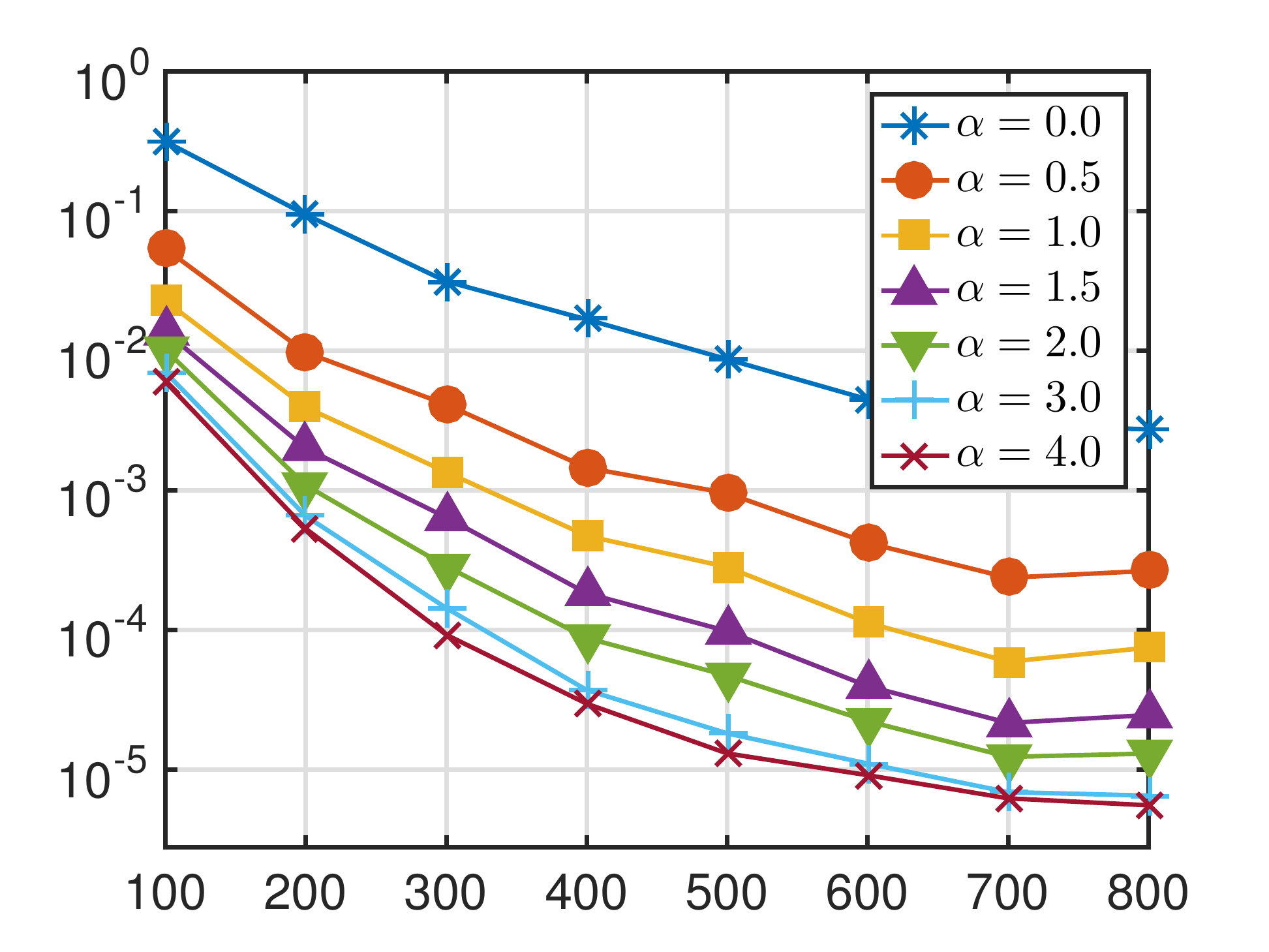} &
\includegraphics[width=5.25cm]{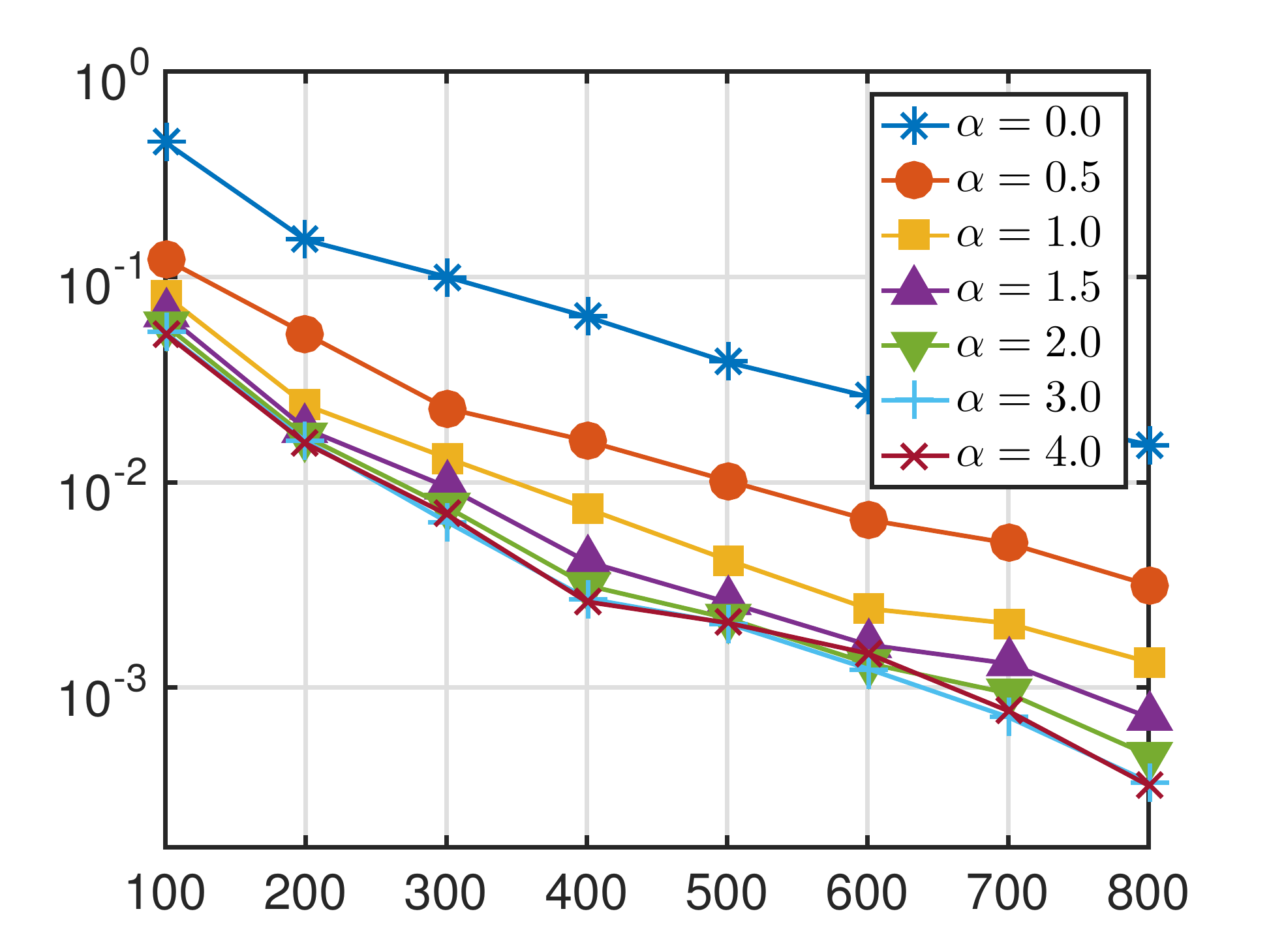} &
\includegraphics[width=5.25cm]{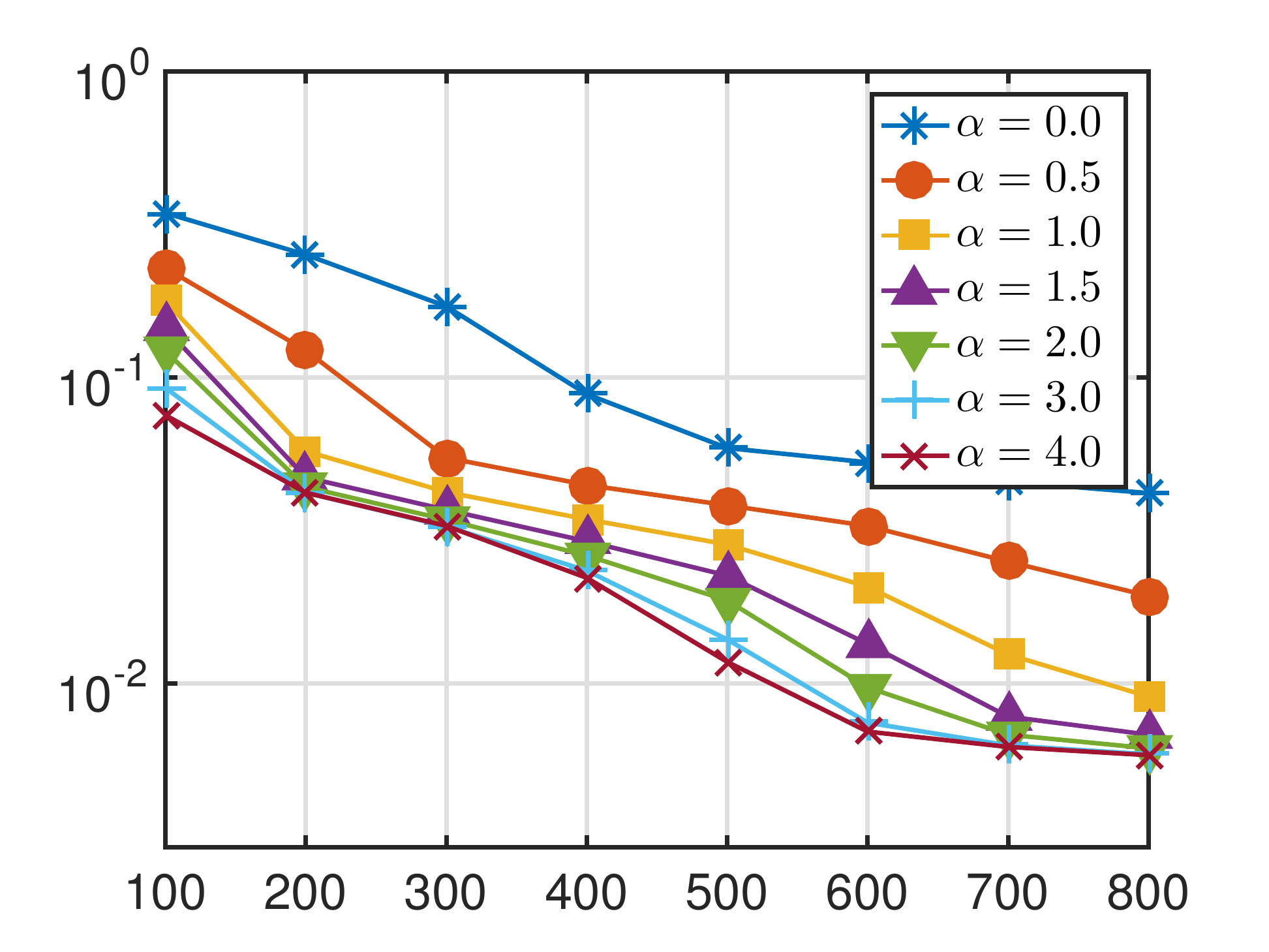}
\\
$d=3$ & $d=5$ & $d=10$
\end{tabular}
\end{center}
\caption{
The error $\| f - \tilde{f} \|_{L^\infty}$ (averaged over $50$ trials) against $m$ for $f(t) = \exp(-(t_1+\ldots+t_d)/d)$ and Legendre polynomials with points drawn from the uniform measure.  A total degree index set of degree $K$ was used, where $(d,K) = (3,24),(5,10),(10,5)$.  The weights were taken to be $w_i = (u_i)^{\alpha}$ for various $\alpha$, where $u_i$ is as in \R{LU_ui}.
}
\label{f:LU_dimcomp}
\end{figure}

Numerical verification of this result is given in Figure \ref{f:LU_dimcomp}.  It is worthwhile noting that a somewhat smaller error in this case can often be be achieved by taking larger weights of the form $w_i = (u_i)^{\alpha}$ with $\alpha >1$.  However, this effect decreases somewhat in higher dimensions.

\vspace{1pc} \noindent
\textit{Tensor Legendre polynomials, random sampling from the Chebyshev measure.}

\cor{
\label{c:LC}
Let $\nu(t) = 2^{-d}$ and $\mu(t) = \prod^{d}_{j=1} \frac{1}{\pi (1-t^2_j)^{1/2}}$.  Then, for any $\Delta \subseteq I_K$ with $| \Delta | \leq s$ we have 
\be{
\label{LC_w1}
\cM(\Delta ; u , 1 ) \leq 2 \times (\pi/2)^{d} (4/\pi)^{\min\{K,d\}} s,
}
provided $I_K = I^{TD}_{K}$ is the total degree space of degree $K$.  Conversely, if $\Delta \subseteq I_K$ is a lower set then
\be{
\label{LC_w2}
 \cM(\Delta ; u , u) \leq 2 \min \left \{ 2^d s , (\pi/2)^{d} s^{\log(1+4/\pi)/\log(2)} \right \}, 
}
regardless of the choice of $I_K$.
}
\prf{
Recall that the weights satisfy \R{u_i_LC_bound} in this case.  The first result follows immediately from this bound.  For the second, we note first that
\bes{
| \Delta |_u \leq (\pi/2)^d K(\Delta),\qquad K(\Delta) = \sum_{i \in \Delta} (4 / \pi)^{|i|_0}.
}
We now claim that $K(\Delta) \leq |\Delta|^{\log(1+4/\pi)/\log(2)}$ for any lower set, thus yielding \R{LC_w2}.  To establish this claim we shall adapt arguments given in \cite{MiglioratiThesis}.  We use induction on $n = | \Delta |$.  If $n=0$ then $\Delta = \{0\}$ (since $\Delta$ is lower) and the claim trivially holds.  Now assume the result holds for $n$ and let $\Delta$ be lower with $|\Delta | = n+1$.   Without loss of generality, $i_1 \neq 0$ for some $i \in \Delta$.  Let $J$ be the maximal value of $i_1$ for $i = (i_1,\ldots,i_d) \in \Delta$ and define the sets
\bes{
 \Delta_k  = \left \{ \hat{i} = (i_2,\ldots,i_d) : (k,i_2,\ldots,i_d) \in \Delta \right \} \subseteq \bbN^{d-1}.
}
Notice each $\Delta_k$ is a lower set and we have the inclusions $\Delta_J \subseteq \Delta_{J-1} \subseteq \cdots \subseteq \Delta_0$.  Since $J \geq 1$ we also have $| \Delta_k | < | \Delta |$ for any $k$, hence the induction hypothesis gives
\be{
\label{mid_step}
K(\Delta) = K(\Delta_0) + 4/\pi \sum^{J}_{k=1} K(\Delta_k) \leq | \Delta_0 |^{\beta} + 4/\pi \sum^{J}_{k=1} | \Delta_k |^{\beta},
}
where $\beta = \log(1+4/\pi)/\log(2) >1$.  We now claim the following.  For $a_0 \geq a_1 \geq \cdots \geq a_n > 0$ it holds that 
\bes{
(a_0+\ldots+a_n)^{\beta} \geq a^{\beta}_0 + 4/\pi (a^{\beta}_1+\ldots+a^{\beta}_n).
}
Proof of this claim follows identical steps to that of \cite[Lem.\ 2.3]{MiglioratiThesis}.  Returning to \R{mid_step}
\bes{
K(\Delta) \leq \left ( \sum^{J}_{k=0} | \Delta_k | \right )^{\beta} = | \Delta |^{\beta},
}
as required.
}

As in the previous examples, this proposition shows that setting the weights $w_i = u_i$ improves the recovery guarantee for lower sets.  To the best of our knowledge, this result has not appeared previously elsewhere.

\rem{
\label{r:logfactors}
Besides the quantity $\cM(\Delta ; u,w)$, the main estimate \R{main_est} also involves the log factor $\log(2N \max \{ \sqrt{| \Delta |_w} , 1 \} )$ (there is also a second log factor depending on the failure probability $\epsilon$, but this is independent of $N$ and $\Delta$ and hence will not discussed further).  In the unweighted case $w_i=1$, since $| \Delta | = s \leq N$ this reduces to a log factor proportional to $\log(2N)$.  Conversely, in the case $w_i = u_i$ one has the log factor
\bes{
\log(2N \sqrt{| \Delta |_u } ) = \log\left( N \sqrt{2\cM(\Delta ; u,u)} \right).
}
Corollaries \ref{c:CC}--\ref{c:LC} can therefore be used to estimate the right-hand side.  In particular, if $\Delta$ is lower, then Corollaries \ref{c:CC} and \ref{c:LU} give a resulting log factor proportional to $\log(2N)$ for the CC and LU cases, since $\cM(\Delta ; u,u)$ is polynomial in $s$ independently of $d$ in these cases, whereas for the LC case Corollary \ref{c:LC} gives a factor proportional to $d + \log(2N)$.
}

\rem{
\label{r:logfactorslower}
When $\Delta$ is lower and the weights are chosen as $w_i = u_i$, the estimates for $\cM(\Delta;u,u)$ in Corollaries \ref{c:CC}--\ref{c:LC} are independent of the choice of truncated space $I_K$ (provided $\Delta \subseteq I_K$).  This choice only affects the parameter $N = |I_K|$, which, as discussed in the previous remark, arises only as a log factor in the measurement condition.  While we have used a total degree space in our numerical experiments for simplicity, a viable alternative (introduced in \cite{ChkifaDownwardsCS})  involves taking $I_K$ to be the union of all lower sets of cardinality $s$.  This is precisely the hyperbolic cross $I^{HC}_{s}$ of order $s$.  Estimates for $N$ in this case are given in \R{HCcardinality}, and lead to a bound for the logarithmic factor of the form $\log(2N) \lesssim \min \left \{ \log(2s) + d , \log(d) \log(2s) \right \}$.
}

\subsection{Support estimation via weighted $\ell^1$ minimization}\label{ss:support_est}
We now turn our attention to a different use of weights: namely, to improve recovery performance when prior information about the support of $x$ is available.  As mentioned, a number of recent works have empirically demonstrated the benefits of this strategy in multivariate polynomial approximation.  In this section, we provide theoretical support to this work.

To do this, we shall assume for simplicity that $u_i = 1$, $\forall i$, although what follows extends to general $u_i$'s.  Let $\Delta \subseteq I$ be the set of coefficients we wish to recover and suppose that $\Gamma \subseteq I$ is an estimate for $\Delta$ based on prior information.  In order to exploit this knowledge, we choose weights
\be{
\label{wi_prior}
w_i = \left \{ \begin{array}{ll} \gamma & i \in \Gamma \\ 1 & i \notin \Gamma \end{array} \right .,
}
where $0 < \gamma < 1$ is a fixed quantity based on the confidence of our estimate.  Define scalars
\be{
\label{rho_sigma_def}
\rho = \frac{| \Delta \cap \Gamma |}{| \Gamma |},\qquad \sigma = \frac{|\Gamma | }{|\Delta |},
}
and observe that $\rho,\sigma \rightarrow 1$ as the accuracy of $\Gamma$ increases.  Then:

\cor{
Let $\Delta,\Gamma \subseteq I$, $| \Delta | =s$, and suppose that $u_i=1$, $\forall i$, and $w = \{w_i \}_{i \in I}$ is as in \R{wi_prior}.  Then, if $\cM$ is as in \R{measurement_condition}, we have
\bes{
\cM(\Delta ; 1,1) = 2 s,
}
and
\bes{
\cM(\Delta \cup \Gamma ; 1, w) = \left ( 2 + \sigma ( 1 + \gamma - 2 \rho ) \right ) s,
}
where $\rho$ and $\sigma$ are as in \R{rho_sigma_def}.  In particular, if 
\bes{
\rho > \frac{1+\gamma}{2},
}
then
\bes{
\cM(\Delta \cup \Gamma ; 1 , w) < \cM(\Delta ; 1,1).
}
}
\prf{
Since $w_i = 1$ for $i \notin \Delta \cup \Gamma$ we have
\eas{
\cM(\Delta \cup \Gamma ; 1,w) &= | \Delta \cup \Gamma | + | \Delta \cup \Gamma |_w
\\
& = | \Gamma | + 2 | \Delta \backslash \Gamma | + \gamma | \Gamma |  
\\
& = | \Gamma | + 2 | \Delta | - 2 | \Delta \cap \Gamma | + \gamma | \Gamma |
\\
& = \left ( \sigma + 2 - 2 \rho \sigma + \gamma \sigma \right ) s
\\
& = \left ( 2 + \sigma ( 1 + \gamma - 2 \rho ) \right ) s,
}
as required.
}

This result implies the following.  If $\gamma$ is sufficiently small and if over half of the support set $\Delta$ is correctly guessed, then the above weighting strategy leads to a smaller measurement condition than in the unweighted case.  In other words, weighting based on sufficiently good prior coefficient estimates can reduce the number of measurements required.

\rem{
Weighted $\ell^1$ minimization with prior support information has been been explored in a number of works \cite{BahWard,FriedlanderWeighted,YuBaekWeighted}.  The setup we consider above is based on that of Friedlander et al.\ \cite{FriedlanderWeighted}.  In most prior works, the measurements are usually taken to be of random Gaussian type, which leads to stronger guarantees than ours.  We are aware of no works that consider prior support information for random sampling of orthonormal systems of functions.  In passing, we note that the improved recovery guarantee of Theorem \ref{t:main} is critical to this analysis, since it allows for arbitrary weights $w_i$.
}

\section{Proof of Theorem \ref{t:main}}\label{s:proof}
We now give the proof of Theorem \ref{t:main}.  As is standard, we first renormalize the problem \R{fin_min} as follows.  Define the infinite matrix $A$
\be{
\label{A_def}
A = \frac{1}{\sqrt{m}} \left \{ \phi_j(t_i) \sqrt{\nu(t_i) / \mu(t_i) } \right \}^{m,\infty}_{i=1,j=1} = \frac{1}{\sqrt{m}} U.
}
Much like $U$, $A$ is a bounded operator from $\ell^1_u(\bbN)$ to $\bbC^m$ whenever the weights are as in \R{u_i_def} (recall Lemma \ref{l:U_bounded}).  Note also that
\bes{
\bbE((A^* A)_{ij}) = \frac{1}{m} \sum^{m}_{k=1} \bbE \left ( \overline{\phi_i(t_k)} \phi_j(t_k) \nu(t_k) / \mu(t_k) \right ) = \int_{D} \overline{\phi_i(t)} \phi_j(t) \nu(t) \D t = \delta_{ij},\quad i=1,\ldots,m,\ j \in I.
}
It follows that \R{fin_min} is equivalent to
\be{
\label{fin_min_A}
\min_{z \in P_K(\ell^1_w(I))} \| z \|_{1,w}\ \mbox{subject to $\| A P_K z - y \| \leq \eta/\sqrt{m}$},
}
where $y = A x + e$, $x$ are the coefficients of $f$ in the basis $\{ \phi_i \}_{i \in I}$ and $e$ is a noise vector satisfying $\| e \| \leq \eta / \sqrt{m}$.  We consider this problem from now on.

The proof now follows a similar route to prievious nonuniform recovery guarantees in CS, although with some significant modifications.  We first show that Theorem \ref{t:main} follows from the existence of a certain dual certificate (Lemma \ref{l:dual_certificate}), and then construct the dual certificate using a variant of the golfing scheme of D.\ Gross \cite{GrossGolfing}.  Technical lemmas required for this construction are presented in \S \ref{ss:tech_lemmas}.  Our argument involves two key novelties.  First, the handling of infinite tails -- an issue that does not arise in most previously-considered (i.e.\ finite-dimensional) CS setups.  Second, the additional complications, and correspondingly refined estimates leading to \R{main_est}, due to the presence of the weights in the optimization problem.

\subsection{Dual certificate}

\lem{
\label{l:dual_certificate}
Let $w = \{ w_i \}_{i \in I}$ be positive weights and $\Delta \subseteq I_K$ be such that $\min_{i \in I_K \backslash \Delta } \{ w_i \} \geq 1$.  Suppose that
\bes{
(i): \| P_{\Delta} A^*  A P_{\Delta} - P_{\Delta} \| \leq \alpha,\qquad (ii):  \max_{i \in I_K \backslash \Delta} \left \{ \|  A e_i \| / w_i \right \} \leq \beta,
}
and that there exists a vector $\rho = W^{-1} P_K A^*  \xi \in P_K(\ell^2(I))$ for some $\xi \in \bbC^m$, where $W = \mathrm{diag}(w_1,w_2,\ldots)$, such that
\bes{
(iii): \| W(P_{\Delta} \rho - \sgn(P_{\Delta} x)) \| \leq \gamma,\quad (iv): \| P^{\perp}_{\Delta} \rho \|_{\infty} \leq \theta,\quad (v): \| \xi \| \leq \lambda \sqrt{| \Delta |_w}, 
}
for constants $0 \leq \alpha , \theta < 1$ and $\beta , \gamma, \lambda \geq 0$ satisfying $\frac{\sqrt{1+\alpha} \beta \gamma}{(1-\alpha)(1-\theta) } < 1$.  Let $x \in \ell^1_w(I)$, $y = A x + e$ with $\| e \| \leq \eta$ and suppose that $\hat{x}$ is a minimizer of the problem
\bes{
\min_{z \in P_K(\ell^1_w(I))} \| z \|_{1,w}\  \mbox{subject to $\| A P_K z - y \| \leq \eta$.}
}
If $\bar{x} \in P_K(\ell^1_w(I))$ is feasible for this problem, i.e.\ $\| A P_K \bar{x} - y \| \leq \eta$, then the estimate
\be{
\label{l1_error_est}
\| \hat{x} - x \| \leq \left ( C_1 + C_2 \lambda\sqrt{| \Delta |_w}\right ) \left ( 2 \eta + \| x - P_K x \|_{1,u} \right ) +  C_2 \left (2 \| x - P_{\Delta} x \|_{1,w} + \| x - \bar{x} \|_{1,w} \right ),
}
holds,
where $C_1 = \left ( 1 + \frac{\gamma}{1-\theta} \right ) C_0$, $C_2 =  \frac{\beta}{1-\theta}\left ( 1 + \frac{\gamma }{1-\theta} \right ) C_0 + \frac{1}{1-\theta}$ and $C_0 = \left ( 1 - \frac{\sqrt{1+\alpha} \beta \gamma }{(1-\alpha)(1-\theta)} \right )^{-1} \frac{\sqrt{1+\alpha}}{1-\alpha}$.
}
\prf{
Let $v = \hat{x} - P_K x \in P_K(\ell^2(I))$.  Then $P_{\Delta} A^* A P_{\Delta} v =  P_{\Delta} A^*  A  v -  P_{\Delta} A^*  A P^{\perp}_{\Delta} v$.  By $(i)$
\bes{
\| (P_{\Delta} A^* A P_{\Delta})^{-1} \| \leq \frac{1}{1-\alpha},
}
and 
\bes{
\|  P_{\Delta} A^*   \|^2 = \|   A P_{\Delta} \|^2 = \|  P_{\Delta} A^*  A P_{\Delta} \| \leq 1 + \alpha.
}
Thus
\eas{
\| P_{\Delta} v \| &\leq \frac{1}{1-\alpha} \|  P_{\Delta} A^*  \| \|  A v \| + \frac{1}{1-\alpha} \|  P_{\Delta} A^*  A   P^{\perp}_{\Delta}  v \|
 \leq \frac{\sqrt{1+\alpha}}{1-\alpha} \left ( \|  A  v \| + \|   A P^{\perp}_{\Delta} v \| \right ).
}
Observe that $\|  A  v \| = \| A \hat{x} - A  P_K x \| \leq 2 \eta + \| A ( x - P_K x ) \|$, and therefore by Lemma \ref{l:U_bounded}
\be{
\label{UN_v}
\|  A  v \| \leq 2 \eta + \| x - P_K x \|_{1,u}.
}
Hence
\bes{
\| P_{\Delta} v \|  \leq \frac{\sqrt{1+\alpha}}{1-\alpha}\left ( 2 \eta + \| x - P_K x \|_{1,u}+ \|   A P^{\perp}_{\Delta} v \| \right ). 
}
The third term can be estimated as follows:
\eas{
\|   A P^{\perp}_{\Delta} v \| \leq \sum_{i \notin \Delta} | v_i | \|   A e_i \| \leq \beta \| P^{\perp}_{\Delta} v \|_{1,w},
}
where the latter inequality is due to $(ii)$.  Hence we get
\be{
\label{v_Delta_Delta_perp}
\| P_{\Delta} v \|  \leq \frac{\sqrt{1+\alpha}}{1-\alpha} \left ( 2 \eta + \| x - P_K x \|_{1,u} +\beta  \|  P^{\perp}_{\Delta} v \|_{1,w} \right ).
}
We shall return to this inequality later, but let us now consider $\hat{x}$.
\ea{
\nm{\hat{x}}_{1,w} &= \nm{P_\Delta\hat{x}}_{1,w} + \nmu{P^{\perp}_\Delta \hat{x} }_{1,w} \nn
\\
& \geq \Re \ip{P_\Delta W \hat{x}}{\mathrm{sign}(P_\Delta x)} + \nmu{P^{\perp}_\Delta  v }_{1,w} - \nmu{P^{\perp}_\Delta x}_{1,w} \nn
\\
&= \Re \ip{P_\Delta W v}{\mathrm{sign}(P_\Delta x)} + \nm{P_\Delta x}_{1,w} + \nmu{P^{\perp}_\Delta v }_{1,w} - \nmu{P^{\perp}_\Delta x}_{1,w} \nn
\\
& = \Re \ip{P_\Delta W v}{\mathrm{sign}(P_\Delta x)} + \nm{x}_{1,w} + \nmu{P^{\perp}_\Delta v }_{1,w}  - 2 \nmu{P^{\perp}_\Delta x}_{1,w}. \label{sign_ineq}
}
Now let $\bar{x} \in \bbC^{K}$ be feasible.  Then $\| \hat{x} \|_{1,w} \leq \| \bar{x} \|_{1,w}$ and we get
\bes{
\|  \bar{x} \|_{1,w} \geq  \Re \ip{P_\Delta W v}{\mathrm{sign}(P_\Delta x)} + \nm{x}_{1,w} + \nmu{P^{\perp}_\Delta v }_{1,w}  - 2 \nmu{P^{\perp}_\Delta x}_{1,w},
}
which after rearranging gives
\be{
\label{2nd_bound}
\| P^{\perp}_\Delta v \|_{1,w} \leq | \ip{P_\Delta W v}{\mathrm{sign}(P_\Delta x)} | + 2 \| P^{\perp}_\Delta x \|_{1,w} + \| x - \bar{x} \|_{1,w}.
}
We next estimate $| \ip{P_\Delta W v}{\mathrm{sign}(P_\Delta x)} |$.  We have
\bes{
| \ip{P_\Delta W v}{\mathrm{sign}(P_\Delta x)} | \leq | \ip{P_\Delta W v}{\mathrm{sign}(P_\Delta x) - P_{\Delta} \rho} | +  | \ip{W v}{\rho} | + | \ip{P^{\perp}_{\Delta} W v}{P^{\perp}_{\Delta} \rho} |.
}
Note that
\bes{
| \ip{P_\Delta W v}{\mathrm{sign}(P_\Delta x) - P_{\Delta} \rho} | \leq \gamma \| P_{\Delta} v \|,
}
and also that
\bes{
\ip{W v}{\rho} = \ip{v}{W \rho} = \ip{v}{A^*  \xi} = \ip{ A v}{\rho}.
}
Hence, \R{UN_v} and $(v)$ give
\bes{
|\ip{W v}{\rho} | \leq \|  A v \| \lambda\sqrt{s} \leq (2 \eta+\| x - P_K x \|_{1,u} ) \lambda\sqrt{|\Delta|_w}.
}
Finally, by $(iv)$, we have
\bes{
| \ip{P^{\perp}_{\Delta} W v}{P^{\perp}_{\Delta} \rho} | \leq \| P^{\perp}_{\Delta} \rho \|_{\infty} \| P^{\perp}_{\Delta} v \|_{1,w} \leq \theta \| P^{\perp}_{\Delta} v \|_{1,w}.
}
Hence 
\bes{
| \ip{P_\Delta W v}{\mathrm{sign}(P_\Delta x)} | \leq \gamma \| P_{\Delta} v \| + (2\eta+ \| x - P_K x\|_{1,u} ) \lambda\sqrt{|\Delta|_w} + \theta \| P^{\perp}_{\Delta} v \|_{1,w},
}
and substituting into \R{2nd_bound} and rearranging yields
\bes{
(1-\theta) \| P^{\perp}_{\Delta} v \|_{1,w} \leq  \gamma \| P_{\Delta} v \| + (2\eta+\| x - P_K x\|_{1,u} ) \lambda\sqrt{|\Delta|_w} + 2 \| P^{\perp}_{\Delta} x \|_{1,w}+ \| x - \bar{x} \|_{1,w}.
}
Applying \R{v_Delta_Delta_perp} now gives
\eas{
\| P_{\Delta} v \| \leq \frac{\sqrt{1+\alpha}}{1-\alpha} \Bigg [&  2 \eta  + \| x - P_K x \|_{1,u}  
\\
& +  \frac{\beta}{1-\theta} \left ( \gamma \| P_{\Delta} v \| +(2\eta+\| x - P_K x\|_{1,u} ) \lambda\sqrt{|\Delta|_w} + 2 \| P^{\perp}_{\Delta} x \|_{1,w}+ \| x - \bar{x} \|_{1,w}  \right ) \Bigg ],
}
and therefore
\eas{
\| P_{\Delta} v \| &\leq \left ( 1 - \frac{\sqrt{1+\alpha} \beta \gamma}{(1-\alpha)(1-\theta)} \right )^{-1}  \frac{\sqrt{1+\alpha}}{1-\alpha} \left ( 1 + \frac{\beta}{1-\theta} \lambda\sqrt{|\Delta|_w} \right ) \left(2\eta+\|x - P_K x\|_{1,u} \right ) 
\\
& +\left ( 1 - \frac{\sqrt{1+\alpha} \beta \gamma }{(1-\alpha)(1-\theta)} \right )^{-1} \frac{\sqrt{1+\alpha} \beta}{(1-\alpha)(1-\theta)} \left ( 2 \| P^{\perp}_{\Delta} x \|_{1,w} + \| x - \bar{x} \|_{1,w} \right )
\\
& = C_0 \left ( 1 + \frac{\beta}{1-\theta} \lambda\sqrt{|\Delta|_w} \right ) \left ( 2\eta + \| x - P_K x \|_{1,u} \right )  + C_0 \frac{\beta}{1-\theta} \left ( 2 \| P^{\perp}_{\Delta} x \|_{1,w} + \| x - \bar{x} \|_{1,w} \right ).
}
Since $w_i \geq 1$, $i \in I_K \backslash \Delta$, we have $\| P^{\perp}_{\Delta} v \| \leq \| P^{\perp}_{\Delta} v \|_1  \leq \| P^{\perp}_{\Delta} v \|_{1,w}$, and hence
\eas{
\| v \| \leq& \| P_{\Delta} v \| +  \| P^{\perp}_{\Delta} v \|_{1,w}
\\
 \leq& \left ( 1 + \frac{\gamma }{1-\theta} \right )  \| P_{\Delta} v \| + \frac{1}{1-\theta} \left ( 2 \eta + \| x - P_K x \|_{1,u} \right ) \lambda\sqrt{|\Delta|_w} + \frac{1}{1-\theta} \left ( 2 \| P^{\perp}_{\Delta} x \|_{1,w} + \| x - \bar{x} \|_{1,w} \right )
\\
\leq& \left [ \left ( 1 + \frac{\gamma }{1-\theta} \right ) C_0 \left (1 + \frac{\beta}{1-\theta} \lambda\sqrt{| \Delta |_w} \right )  + \frac{1 }{1-\theta} \lambda \sqrt{|\Delta|_w} \right ] \left (2 \eta+ \| x - P_K x \|_{1,u} \right )  
\\
&+ \left [ \frac{C_0 \beta}{1-\theta}  \left ( 1 + \frac{\gamma}{1-\theta} \right ) + \frac{1}{1-\theta} \right ]  \left ( 2 \| P^{\perp}_{\Delta} x \|_{1,w} + \| x - \bar{x} \|_{1,w} \right ),
}
as required.
}

\subsection{Technical lemmas}\label{ss:tech_lemmas}
For the dual certificate construction in \S \ref{ss:dual_constr} we first require the following four technical lemmas:

\lem{
\label{t:near_isometry}
Let $\Delta \subseteq I$, $1 \leq | \Delta | < \infty$, $0 < \epsilon < 1$, $\delta>0$ and $u = \{ u_i \}_{i \in I}$ be as in \R{u_i_def}.  Then
\bes{
\| P_{\Delta} A^*  A P_{\Delta} - P_{\Delta} \| \leq \delta ,
}
with probability at least $1-\epsilon$, provided
\bes{
m \geq  | \Delta |_u \cdot  \log \left ( 2 |\Delta| / \epsilon \right ) \cdot \left ( 2 \delta^{-2} + 2 \delta^{-1}/3 \right ).
}
}
\prf{
Let $Y_i = \{ \sqrt{\nu(t_i)/\mu(t_i)} \phi_j(t_i) \}_{j \in \Delta}$ and observe that
\bes{
P_{\Delta} A^* A P_{\Delta} - P_{\Delta} = \sum^{m}_{i=1} X_i,
}
where $X_i = \frac1m \left ( Y_i Y^*_i - I  \right ) \in \bbC^{|\Delta| \times |\Delta|}$ satisfies $\bbE(X_i) = 0$.  Note that
\bes{
\| Y_i \|^2 = \sum_{j \in \Delta} \frac{\nu(t_i)}{\mu(t_i)} | \phi_j(t_i) |^2 \leq \sum_{j \in \Delta} u^2_j = | \Delta |_u.
}
Hence
\bes{
\| X_i \| = \sup_{\| x \|=1} | \ip{X_i x}{x} | = m^{-1} \sup_{\| x \|=1} \left | |\ip{Y_i}{x} |^2 - 1 \right | \leq m^{-1} |\Delta |_u ,
}
where in the last step we use the fact that $| \Delta |_u \geq | \Delta | \geq 1$.  Also
\bes{
\bbE (X^2_i) = m^{-2} \bbE \left ( (\| Y_i \|^2 - 2 )Y_i  Y^*_i + I \right ),
}
Since $\bbE(Y_i Y^*_i) = I$ we have
\bes{
|\ip{\bbE (X^2_i) x}{x} |\leq m^{-2} | \Delta |_u  \| x \|^2,
}
and therefore
\bes{
\nm{\sum^{m}_{i=1} \bbE (X^2_i) }=  \sup_{\| x \|=1}\left | \sum^{m}_{i=1} \ip{\bbE(X^2_i) x}{x} \right | \leq m^{-1}  | \Delta |_u .
}
The result now follows immediately from the matrix Bernstein inequality \cite[Cor.\ 8.15]{FoucartRauhutCSbook}.  
}

\lem{
\label{l:near_isometry_vector}
Let $\Delta \subseteq I$, $1 \leq | \Delta | < \infty$, $0 < \epsilon < \E^{-1}$, $\delta>0$, $u = \{ u_i \}_{i \in I}$ be as in \R{u_i_def} and $z \in \ell^2(I)$.  Then
\bes{
\| (P_{\Delta} A^*  A P_{\Delta} - P_{\Delta}) z \| \leq \delta \| z \| ,
}
with probability at least $1-\epsilon$, provided
\bes{
m \geq | \Delta |_u  \cdot  \log(\epsilon^{-1})  \cdot \left ( 8 \delta^{-2} + 28 \delta^{-1} / 3 \right ).
}
}
\prf{
Let $\| z \|=1$ without loss of generality.  As in the previous proof, write $P_{\Delta} A^*  A P_{\Delta} - P_{\Delta} = m^{-1} \sum^{m}_{i=1} ( Y_i Y^*_i - I )$ so that
\bes{
\nm{(P_{\Delta} A^*  A P_{\Delta} - P_{\Delta}) z} = \nm{\sum^m_{i=1} Z_i },
}
where $Z_i = m^{-1} (Y_i Y^*_i - I ) z$.  Observe that $\bbE(Z_i) = 0$ and the $Z_i$ are independent copies of a single random vector.  Arguing as in the previous lemma, we note that
\bes{
\| Z_i \| \leq m^{-1} | \Delta |_u,
}
and
\bes{
\bbE \| Z_i \|^2 \leq m^{-2} \max\{ | \Delta |_u,1\},\qquad \sup_{\| x \|=1} \bbE | \ip{Z_i}{x} |^2 \leq m^{-2} | \Delta |_u.
}
Suppose that $m \geq 4 | \Delta |_u \delta^{-2}$ so that $\sqrt{m \bbE \| Z_i \|^2} \leq \delta / 2$.  It now follows from \cite[Cor.\ 8.45]{FoucartRauhutCSbook} that
\bes{
\bbP \left ( \nm{(P_{\Delta} A^*  A P_{\Delta} - P_{\Delta}) z} > \delta \right ) \leq \exp \left ( -\frac{m}{\max \{ | \Delta |_u,1\}} \frac{\delta^2/8}{1+ 7 \delta / 6 } \right ),
}
which gives the result.
}

\lem{
\label{l:second_term_est}
Let $0 < \epsilon < 1$, $\delta > 0$, $w = \{ w_i \}_{i \in I}$ be weights, $\Delta \subseteq I_K$ and suppose that $\min_{i \in I_K \backslash \Delta } \{ w_i \} \geq 1$.  Then
\be{
\label{second_term_est_m}
\max_{i \in I_K \backslash \Delta } \left \{ \| A e_i \| / w_{i} \right \} \leq \sqrt{1+\delta},
}
with probability at least $1-\epsilon$, provided
\bes{
m \geq 2 \max_{i \in I_K \backslash \Delta} \{ u^2_i / w^2_i \}  \cdot  \log\left ( 2 N / \epsilon \right ) \cdot \left ( \delta^{-2} + \delta^{-1}/3 \right ),
}
where $N = |I_K|$ and $u = \{u_i\}_{i \in I}$ is as in \R{u_i_def}.
}
\prf{
Fix $i \in I_K \backslash \Delta$.  Then
\eas{
\| A e_i \|^2/w^2_i &= e^*_i A^* A e_i /w^2_i = \frac{1}{m} \sum^{m}_{j=1} \frac{\nu(t_j)}{\mu(t_j)} | \phi_i(t_j) |^2 / w^2_i  \leq \left | \sum^{m}_{j=1} X_j \right | + 1/w^2_i \leq \left | \sum^{m}_{j=1} X_j \right | + 1,
}
where $X_j = m^{-1} \left ( \frac{\nu(t_j)}{\mu(t_j)} | \phi_i(t_j) |^2 - 1 \right )/w^2_i$.  Note that $\bbE(X_j) = 0$ and the $X_j$'s are independent.  Moreover,
\bes{
|X_j| \leq m^{-1} u^2_i / w^2_i ,
}
and
\eas{
\sum^{m}_{j=1}\bbE(|X_j|^2) &= \frac{1}{m w^4_i} \int_{D} \left ( \frac{\nu(t)}{\mu(t)} | \phi_i(t) |^2 - 1 \right )^2 \mu(t) \D t 
\\
&=  \frac{1}{m w^4_i} \left ( \int_{D} \frac{\nu(t)^2}{\mu(t)^2} | \phi_i(t) |^4 \mu(t) \D t  - 1 \right )
\\
& \leq  \frac{u^2_i  }{m w^2_i}.
}
Hence, by Bernstein's inequality,
\bes{
\bbP(\| A e_i \| / w_i \geq \sqrt{1+\delta} ) \leq 2 \exp \left ( -\frac{1}{\kappa} \frac{\delta^2/2}{1+\delta/3} \right ),\quad \forall i \in I_K \backslash \Delta,
}
where $\kappa = m^{-1} \max_{i \in I_R \backslash \Delta } \{ u^2_i / w^2_i \}$.  The result now follows from the union bound.
}

\lem{
\label{l:Sloan}
Let $0 < \epsilon < 1$, $\delta > 0$,  $z \in \ell^2(I)$, $w = \{ w_i \}_{i \in I}$ be weights and $\Delta \subseteq I_K$.  Then
\bes{
\| P_K P^{\perp}_{\Delta} W^{-1} A^* A P_{\Delta} z \|_{\infty} \leq \delta \| z \|.
}
with probability at least $1-\epsilon$, provided
\be{
\label{l:Sloan_m}
m \geq 2 \left (\max_{i \in I_K \backslash \Delta} \{ u^2_i / w^2_i \} \delta^{-2} +\sqrt{|\Delta|_u} \max_{i \in I_K \backslash \Delta} \{ u_i / w_i \} \delta^{-1} /3 \right ) \cdot \log(2 N / \epsilon),
}
where $N = |I_K|$ and $u = \{u_i\}_{i \in I}$ is as in \R{u_i_def}.
}
\prf{
Let $\| z \|=1$ without loss of generality.  Observe that 
\bes{
\| P_K P^{\perp}_{\Delta} W^{-1} A^* A P_{\Delta} z \|_{\infty} = \max_{i \in I_K \backslash \Delta} \frac{|\ip{e_i}{A^* A P_{\Delta} z } |}{w_i}.
}
Fix $i \in I_K \backslash \Delta$.  Then
\bes{
\frac{|\ip{e_i}{A^* A P_{\Delta} z } |}{w_i} = \left | \sum^{m}_{j=1} X_j \right |,
}
where $X_j$ is the random variable
\bes{
X_j = \frac{1}{w_i m} \phi_i(t_j) \left ( \sum_{k \in \Delta} \phi_k(t_j) z_k \right ) \frac{\nu(t_j)}{\mu(t_j)}.
}
Observe that  the $X_j$ are independent and $\bbE(X_j) = 0$ since $i \notin \Delta$.  Also
\bes{
| X_j | \leq \frac{u_i}{w_i m} \left | \sum_{k \in \Delta} \sqrt{\frac{\nu(t_j)}{\mu(t_j)}} \phi_k(t_j) z_k \right | \leq \frac{u_i}{w_i m} \sqrt{|\Delta|_u},
}
and
\bes{
\sum^{m}_{j=1} \bbE(|X_j|^2) \leq \left ( \frac{u_i}{w_i m} \right )^2 m \int_{D} \left | \sum_{k \in \Delta} \phi_k(t) z_k \right |^2 \nu(t) \D t = \frac{u^2_i}{m w^2_i} .
}
Therefore, by Bernstein's inequality and the union bound,
\bes{
\bbP \left ( \| P_K P^{\perp}_{\Delta} W^{-1} A^* A P_{\Delta} z \|_{\infty} > \delta \right ) \leq 2 N \exp \left ( - \frac{m \delta^2/2}{\max_{i \in I_K \backslash \Delta} \{ u^2_i / w^2_i \} + \delta \sqrt{|\Delta|_u} \max_{i \in I_K \backslash \Delta} \{ u_i / w_i \} / 3} \right ).
}
The result now follows immediately.
}

\subsection{Construction of the dual certificate $\rho$}\label{ss:dual_constr}
We are now ready to construct a dual certificate $\rho$ satisfying the conditions of Lemma \ref{l:dual_certificate}.  For parameters, we choose the following values:
\be{
\label{param_choices}
\alpha = 1/4,\quad \beta = \sqrt{5/4},\quad \gamma = 1/8,\quad \theta = 1/2.
}

\vspace{1pc} \noindent \textit{Setup.}
Let $L \in \bbN$ and suppose that $m_1,\ldots,m_L$ are such that $m_1+\ldots+m_L = m$.  If $U$ is as in \R{U_def}, then write $U^{(1)} \in \bbC^{m_1 \times \infty}$ for the submatrix of the first $m_1$ rows, $U^{(2)} \in \bbC^{m_2 \times \infty}$ for the submatrix of the first $m_2$ rows, and so on.  Define $\rho^{(0)} = 0$,
\bes{
\rho^{(l)} =  m^{-1}_l W^{-1} P_K (U^{(l)})^* U^{(l)} P_{\Delta} W\left ( \sgn(P_{\Delta}(x)) - P_{\Delta} \rho^{(l-1)} \right ) + \rho^{(l-1)},\quad l=1,2.
}
and
\bes{
v^{(l)} = W \left ( \sgn(P_{\Delta}(x)) - P_{\Delta} \rho^{(l)} \right ),\quad l=0,1,2.
}
Let $\Theta^{(1)} = \{1\}$ and $\Theta^{(2)} = \{1,2\}$.  For $l=3,\ldots,L$, we define $\Theta^{(l)}$ as follows:
\bull{
\item If 
\eas{
\| (P_{\Delta} - m^{-1}_{l} P_{\Delta} (U^{(l)})^* U^{(l)} P_{\Delta}) v^{(l-1)} \| &\leq a_l \| v^{(l-1)} \|
\\
\| m^{-1}_{l} P_K P^{\perp}_{\Delta} W^{-1} (U^{(l)})^* U^{(l)} P_{\Delta} v^{(l-1)} \|_{\infty} &\leq b_l \| v^{(l-1)} \|
}
for constants $a_l$ and $b_l$ then set $\Theta^{(l)} = \Theta^{(l-1)} \cup \{ l \}$ and
\bes{
\rho^{(l)} = m^{-1}_{l} W^{-1} P_K (U^{(l)})^* U^{(l)} v^{(l-1)} + \rho^{(l-1)},\qquad v^{(l)} = W \left ( \sgn(P_{\Delta}x) - P_{\Delta} \rho^{(l)} \right ).
}
\item Otherwise, set $\Theta^{(l)} = \Theta^{(l-1)}$, $\rho^{(l)} = \rho^{(l-1)}$ and $v^{(l)} = v^{(l-1)}$.
}
We now define the events $A_1,A_2,B_1,B_2,C,D$ as follows:
\eas{
&A_l :\quad \| (P_{\Delta} - m^{-1}_{l} P_{\Delta} (U^{(l)})^* U^{(l)} P_{\Delta}) v^{(l-1)} \| \leq a_l \| v^{(l-1)} \|,\quad l=1,2,
\\
& B_l : \quad \| m^{-1}_{l} P_K P^{\perp}_{\Delta} W^{-1} (U^{(l)})^* U^{(l)} P_{\Delta} v^{(l-1)} \|_{\infty} \leq b_l \| v^{(l-1)} \|,\quad l=1,2,
\\
& C : \quad | \Theta^{(L)} | \geq R
\\
& D : \quad\| P_{\Delta} A^* A P_{\Delta} - P_{\Delta} \| \leq 1/4,
\\
& E : \quad \sup_{i \notin \Delta} \{ \| A e_i \|/w_i \} \leq \sqrt{5/4},
\\
& F : \quad A_1 \cap A_2 \cap B_1 \cap B_2 \cap C \cap D \cap E,
}
where $| \Theta^{(L)} |$ is the cardinality of $\Theta^{(L)}$.  If event $F$ occurs, write $\Theta^{(L)} = \left \{ \tau(1) , \tau(2) , \ldots, \tau(R) , \ldots \right \}$, where the function $\tau$ satisfies $\tau(l) \geq l$ for all $l$, and define the dual certificate as $\rho = \rho^{(\tau(R))}$.

\vspace{1pc} \noindent \textit{Choice of the parameters.}
The idea of the proof is to choose $L$, $R$, $ m_1,\ldots,m_L$, $a_1,\ldots,a_L$ and $b_1,\ldots,b_L$ so that conditions (i)--(v) of Lemma \ref{l:dual_certificate} are fulfilled for the parameter choices \R{param_choices}.  We make the following choices for these parameters.  Write $s = | \Delta |_w$, $s^* = \max \{ s , 1 \}$ and set
\be{
\label{R_choice}
R = \lceil \log_2(8 N \sqrt{s^*}) \rceil,
}
where $N = |I_K|$,
\be{
\label{L_choice}
L = 2+ \lceil \log(7\epsilon^{-1}) \rceil + 10 R,
}
\be{
\label{a_choice}
a_1 = a_2 = \frac{1}{2 \sqrt{\log_2(8 N \sqrt{s^*})}},\quad a_l = 1/2,\quad l=3,\ldots,L.
}
\be{
\label{b_choice}
b_1 = b_2 = \frac{1}{4 \sqrt{s}},\quad b_l = \frac{\log_2(8 N \sqrt{s^*})}{4 \sqrt{s}} ,\quad l=3,\ldots,L.
}
and
\bes{
m_1 = m_2 = \frac{m}{4},\quad m_l = \frac{m}{2(L-2)},\quad l=3,\ldots,L.
}

\vspace{1pc} \noindent \textit{Claim: Event $D$ implies conditions (i)--(v).}
Suppose that event $D$ occurs.  Immediately, events $D$ and $E$ give that conditions (i) and (ii) hold with $\alpha = 1/4$ and $\beta = \sqrt{5/4}$.  Now consider condition (iii).  If $\tau(k-1),\tau(k) \in \Theta^{(L)}$ then note that
\ea{
v^{(\tau(k))} &= W \left ( \sgn(P_{\Delta}(x)) - P_{\Delta} m^{-1}_{\tau(k)} W^{-1} (U^{(\tau(k))})^* U^{(\tau(k))} v^{\tau(k-1)} - P_{\Delta} \rho^{(\tau(k-1))} \right ) \nn
\\
& = \left ( P_{\Delta} - m^{-1}_{\tau(k)} P_{\Delta} (U^{(\tau(k))})^* U^{(\tau(k))} P_{\Delta} \right ) v^{(\tau(k-1))}. \label{v_relation}
}
Hence
\be{
\label{v_bound}
\| v^{(\tau(k))} \| \leq a_{\tau(k)} \| v^{(\tau(k-1))} \| \leq \| v^{(0)} \| \prod^{k}_{j=1} a_{\tau(j)}  \leq \sqrt{s} \prod^{k}_{j=1} a_{\tau(j)}.
}
Observe that
\be{
\label{a_prod}
\prod^{k}_{j=1} a_{\tau(j)} = \frac{1}{2^k \log_2(8 N \sqrt{s^*})} \leq \frac{1}{2^k} .
}
Hence setting $k = R$ in \R{v_bound} and noticing that
\bes{
W (P_{\Delta} u - \sgn(P_{\Delta}x) ) = W  (P_{\Delta} \rho^{(\tau(R))} - \sgn(P_{\Delta}x) ) = v^{(\tau(R))},
}
gives that
\be{
\label{Propcond3}
\| W (P_{\Delta} u - \sgn(P_{\Delta}x) ) \| \leq \sqrt{s} \prod^{R}_{j=1} a_{\tau(j)} = \frac{\sqrt{s}}{2^R} \leq \frac{1}{8}.
}
Thus condition (iii) holds with $\gamma = 1/8$ as required.  Now consider condition (iv).  Observe that
\bes{
P^{\perp}_{\Delta} \rho^{(\tau(k))} = m^{-1}_{\tau(k)} P_K P^{\perp}_{\Delta} W^{-1} (U^{(\tau(k))})^* U^{(\tau(k))} v^{(\tau(k-1))} + P^{\perp}_{\Delta} \rho^{(\tau(k-1))}.
}
Therefore
\bes{
\| P^{\perp}_{\Delta} \rho^{(\tau(k))} \|_{\infty} \leq b_{\tau(k)} \| v^{(\tau(k-1))} \| + \| P^{\perp}_{\Delta} \rho^{(\tau(k-1))} \|_{\infty} \leq \sqrt{s} b_{\tau(k)} \prod^{k-1}_{j=1} a_{\tau(j)} + \| P^{\perp}_{\Delta} \rho^{(\tau(k-1))} \|_{\infty},
}
where we use the convention that $\prod^{k-1}_{j=1} a_{\tau(j)} = 1$ when $k=1$.  Hence
\be{
\label{Propcond4}
\| P^{\perp}_{\Delta} u \|_{\infty} \leq \sqrt{s} \sum^{R}_{k=1} b_{\tau(k)} \prod^{k-1}_{j=1} a_{\tau(j)}.
}
Substituting the values of $a_l$ and $b_l$ into the right-hand side of \R{Propcond4} and using \R{a_prod} gives
\bes{
\sqrt{s} \sum^{R}_{k=1} b_{\tau(k)} \prod^{k-1}_{j=1} a_{\tau(j)} \leq \frac{1}{4} \left ( 1 + \frac{1}{2} + \frac{1}{4 } +  \frac{1}{8} + \ldots + \frac{1}{2^{R-1}}  \right ) \leq \frac{1}{2}.
}
Hence condition (iv) holds with $\theta = 1/2$, as required.

Finally consider condition (v).  Write $\rho^{(\tau(k))} = W^{-1} P_K A^* \xi^{(\tau(k))}$, where
\bes{
\xi^{(\tau(k))} = \frac{\sqrt{m}}{m_{\tau(k)}} U^{(\tau(k))} v^{(\tau(k-1))} + \xi^{(\tau(k-1))}.
}
It follows that
\be{
\label{wk_bound}
\| \xi^{(\tau(k))} \| \leq \frac{\sqrt{m}}{m_{\tau(k)}} \| U^{(\tau(k))} v^{\tau(k-1)} \| + \| \xi^{(\tau(k-1))} \|.
}
Consider the first term on the right-hand side.  We have
\eas{
m^{-1}_{\tau(k)}\| U^{(\tau(k))}  v^{\tau(k-1)} \|^2 &= m^{-1}_{\tau(k)} \ip{P_{\Delta} (U^{(\tau(k))})^* U^{(\tau(k))}  P_{\Delta} v^{(\tau(k-1))}}{v^{(\tau(k-1))}}
\\
& = \| v^{(\tau(k-1))} \|^2 - \ip{v^{(\tau(k))}}{v^{(\tau(k-1))}}
\\
& \leq \| v^{(\tau(k-1))} \|^2 + \| v^{(\tau(k-1))} \| \| v^{(\tau(k))} \|,
}
where in the middle step we use \R{v_relation}.  By \R{v_bound}, it now follows that
\bes{
m^{-1}_{\tau(k)}\| U^{(\tau(k))}  v^{\tau(k-1)} \|^2 \leq s \left ( a_{\tau(k)} + 1 \right ) \left ( \prod^{k-1}_{j=1} a_{\tau(j)} \right )^2,
}
and therefore, returning to \R{wk_bound} and summing over $k=1,\ldots,R$, we get 
\be{
\label{Propcond5}
\| \xi \| \leq  \sqrt{s} \sum^{R}_{k=1} \sqrt{\frac{m}{m_{\tau(k)}}}\sqrt{a_{\tau(k)} + 1} \prod^{k-1}_{j=1} a_{\tau(j)},
}
where $\xi = \xi^{(\tau(R))}$ is such that $u = W^{-1} P_K A^* \xi$.  Now notice that
\bes{
\sqrt{\frac{m}{m_{\tau(k)}}} \sqrt{a_{\tau(k)}+1} \leq \sqrt{6},\quad k=1,2,
}
and
\bes{
\sqrt{\frac{m}{m_{\tau(k)}}} \sqrt{a_{\tau(k)}+1} \leq \sqrt{3L} \qquad k=3,\ldots,R.
}
It follows from \R{Propcond5} and \R{a_prod} that
\bes{
\| w \| \leq \sqrt{s} \left ( \sqrt{6}(1+1/2) + \sqrt{3L} \sum^{R}_{k=3} \frac{1}{2^k \log_2(8 N \sqrt{s^*})} \right ),
}
and therefore condition (v) holds with
\be{
\label{lambda_const_value}
\lambda \lesssim 1 + \frac{\sqrt{\log(\gamma^{-1})}}{\log(2N\sqrt{\max \{ | \Delta |_w , 1 \}})}.
}

\subsection{Event $F$ holds with high probability}
We now derive conditions on $m$ for event $F$ to hold with probability at least $1 -\epsilon$, where $0 < \epsilon < \E^{-1}$.  By the union bound, it suffices to prove that events $A_1,A_2,B_1,B_2,C,D,E$ occur with probability at least $1-\gamma$, where $\gamma = \epsilon / 7$.

\vspace{1pc} \noindent \textit{Events $A_1,A_2$ and $B_1,B_2$.}
We apply Lemma \ref{l:near_isometry_vector} with $\epsilon = \gamma$, $\delta = 1/(2 \sqrt{\log_2(8 N \sqrt{s^*})})$ and $m = m_l = m/4$.  This gives that $\bbP(A^c_l) \leq \gamma$ for $l=1,2$ provided $m$ satisfies
\be{
\label{A12_mcond}
m \gtrsim | \Delta |_u \cdot \log(\gamma^{-1}) \cdot \log_2(8 N \sqrt{s^*}).
}
Similarly, for events $B_1$ and $B_2$ we apply Lemma \ref{l:Sloan} with $\epsilon = \gamma$, $\delta = 1/(4 \sqrt{s})$ and $m = m_l = m/4$.  This gives that $\bbP(B^c_l) \leq \gamma$ for $l=1,2$, provided
\bes{
m \gtrsim \left (  |\Delta|_w \max_{i \in I_K\backslash \Delta} \{ u^2_i / w^2_i \} + \sqrt{|\Delta|_u |\Delta|_w} \max_{i \in I_K\backslash \Delta} \{ u_i / w_i \} \right ) \cdot \log(2N / \gamma).
}
After simplifying, we see that it suffices that
\be{
\label{B12_mcond}
m \gtrsim  \left ( | \Delta |_u + \max_{i \in I_K \backslash \Delta} \{ u^2_i / w^2_i \} |\Delta|_w \right ) \cdot \log(2N/\gamma). 
}

\vspace{1pc} \noindent \textit{Event $C$.}
Define the random variables $X_{1},\ldots,X_{L-2}$ by 
\bes{
X_{l} = \left \{ \begin{array}{ll} 1 & v^{(l+2)} \neq v^{(l+1)} \\ 0 & \mbox{otherwise} \end{array} \right . ,
}
so that $\bbP(C^{c}) = \bbP(| \Theta^{(L)} | < R ) = \bbP(X_{1}+\ldots+X_{L-2} < R )$.\footnote{In some of the first presentations of the golfing scheme \cite{Candes_Plan,GrossGolfing}, it was assumed that these random variables were independent, which is not the case in general.  This issue was fixed in \cite{BAACHGSCS} via a more careful argument.  Here we follow the approach of \cite{KrahmerGrossKuengPhaseLift}.}  Observe that
\bes{
 \bbP\left (\sum^{L-2}_{l=1}X_{l}< R \right ) = \bbE \left ( \bbP \left ( X_{L-2} < R - \sum^{L-3}_{l=1}X_{l} | X_{L-3},\ldots,X_1 \right ) \right ).
}
When conditioned on an instance of $X_1,\ldots,X_{L-3}$ the variable $X_{L-2}$ has a Bernoulli distribution with some parameter $p(X_1,\ldots,X_{L-3})$.  If $X$ is a Bernoulli random variable with parameter $p$, then the function $\bbP(X<t)$ is a nonincreasing function of $p$ for any fixed $t \in \bbR$.  It follows that 
\be{
\label{reduction_to_independent}
 \bbP\left (\sum^{L-2}_{l=1}X_{l}< R \right ) \leq \bbP \left ( X'_{L-2} + \sum^{L-3}_{l=1} X_l < R \right ),
}
where $X'_{L-2}$ is an independent Bernoulli random variable with parameter $p'$ satisfying
\bes{
p' \leq \min_{x_1,\ldots,x_{L-3} \in \{0,1\}} p(x_1,\ldots,x_{L-3}).
}
We now wish to guarantee that $p' \geq 9/10$.  Observe that $X_{L-2} = 0$ if, for $l=L$, either of the following events occur:
\eas{
C_1 &: \| (P_{\Delta} - m^{-1}_{l} P_{\Delta} (U^{(l)})^* U^{(l)} P_{\Delta}) v^{(l-1)} \| > a_l \| v^{(l-1)} \|
\\
C_2 &: \| m^{-1}_{l} P_K P^{\perp}_{\Delta} W^{-1} (U^{(l)})^* U^{(l)} P_{\Delta} v^{(l-1)} \|_{\infty} > b_l \| v^{(l-1)} \|.
}
Applying Lemma \ref{l:near_isometry_vector} with $m=m_L$, $\delta = 1/2$, $\epsilon = 1/20$ and Lemma \ref{l:Sloan} with $m=m_L$, $\delta = \log_2(8 N \sqrt{s^*}) / (4 \sqrt{s})$ and $\epsilon = 1/20$ we now see that $p' \geq 9/10$, provided $m$ satisfies
\bes{
m_L \gtrsim | \Delta |_u ,
}
and
\bes{
m_L \gtrsim \left (  |\Delta|_w \max_{i \in I_K\backslash \Delta} \{ u^2_i / w^2_i \} + \sqrt{|\Delta|_u |\Delta|_w} \max_{i \in I_K\backslash \Delta} \{ u_i / w_i \}  \right ) \log(40K)) / \log_2(8 N \sqrt{s^*}).
}
Since $m_L = m/(2(L-2))$, these now reduce to
\be{
\label{E1_cond}
m \gtrsim | \Delta |_u \cdot \left ( \log(\epsilon^{-1}) + \log_2(8 N \sqrt{s^*}) \right ),
}
and
\be{
\label{E2_cond}
m \gtrsim \left ( | \Delta |_u + \max_{i \in I_K \backslash \Delta} \{ u^2_i / w^2_i \} | \Delta |_w \right ) \cdot \log(2N) \cdot \log(\epsilon^{-1}),
}
respectively.  

Assuming now that \R{E1_cond} and \R{E2_cond} hold, we have $p' \geq 9/10$.  Moreover, repeating the same arguments we can iterate the estimate \R{reduction_to_independent} to obtain
\bes{
 \bbP\left (\sum^{L-2}_{l=1}X_{l}< R \right ) \leq \bbP \left ( \sum^{L-2}_{l=1} X'_l < R \right ),
}
where the $X'_l$ are independent Bernoulli random variables with parameter $9/10$.  Following the arguments of \cite{BAACHGSCS}, we deduce that
\bes{
 \bbP\left (\sum^{L-2}_{l=1}X_{l}< R \right ) \leq \exp \left ( -2 (L-2) \left ( \frac{9}{10} - \frac{R}{L-2} \right )^2 \right ).
}
By construction
\bes{
2 (L-2) \left ( \frac{9}{10} - \frac{R}{L-2} \right )^2 \geq 2 \log(\gamma^{-1}) \left ( \frac{9}{10} - \frac{1}{10} \right )^2 > \log(\gamma^{-1}).
}
Hence 
\bes{
\bbP(C^c) =  \bbP\left (\sum^{L-2}_{l=1}X_{l}< R \right ) \leq \gamma,
}
as required.

\vspace{1pc} \noindent \textit{Events $D$ and $E$.}
For event $D$, we apply Lemma \ref{t:near_isometry} with $\epsilon = \gamma$ and $\delta = 1/4$.  This yields that if
\be{
\label{D_cond}
m \gtrsim | \Delta |_u \cdot \log(2 | \Delta | / \gamma ),
}
then $\bbP(D^c) \leq \gamma$.  Similarly, suppose that
\be{
\label{E_cond}
m \gtrsim \max_{i \in I_R \backslash \Delta } \{ u^2_i / w^2_i \} \cdot \log(2N/\gamma),
}
then Lemma \ref{l:second_term_est} gives that $\bbP(E^c) \leq \gamma$.

\subsection{Proof of Theorem \ref{t:main}}
Recall that \R{fin_min} is equivalent to \R{fin_min_A}.  The estimate now follows from Lemma \ref{l:dual_certificate}, provided conditions (i)--(v) hold.  As demonstrated in the previous section, (i)--(v) hold with probability at least $1-\epsilon$ and with $\lambda$ given by \R{lambda_const_value}, provided \R{A12_mcond}, \R{B12_mcond}, \R{E1_cond}, \R{E2_cond}, \R{D_cond} and \R{E_cond} hold.  However, due to the assumption on $\epsilon$ and the choice of $\gamma$, these are all implied by \R{main_est}.

\section*{Acknowledgements}
The work was supported in part by the Natural Sciences and Engineering Research Council of Canada through grant 611675 and an Alfred P. Sloan Research Fellowship.  The author would particularly like to thank Abdellah Chkifa, Clayton Webster, Hoang Tran and Guannan Zhang for introducing him to the concept of lower sets.  The results of \S \ref{ss:recovery_examples} are due to their insight.  He would also like to like to thank Rick Archibald, Nilima Nigam, Clarice Poon and Tao Zhou for useful discussions and comments.

\bibliographystyle{abbrv}
\small
\bibliography{CSFunInterpRefs}

\end{document}